\documentclass{article}

\usepackage{arxiv}

\usepackage[utf8]{inputenc} 
\usepackage[T1]{fontenc}    
\usepackage{hyperref}       
\usepackage{url}            
\usepackage{booktabs}       
\usepackage{amsfonts}       
\usepackage{nicefrac}       
\usepackage{microtype}      
\usepackage{lipsum}
\usepackage{mathtools}
\usepackage{amsmath}
\usepackage{graphicx}
\usepackage{algpseudocode}
\usepackage{varwidth}
\usepackage{floatrow}
\newfloatcommand{capbtabbox}{table}[][\FBwidth]
\usepackage{caption}
\usepackage{subcaption}

\title{Some Numerical Simulations Based on Dacorogna Example Functions in Favor of Morrey Conjecture}

\author{
 Xinghao Dong \\
  Department of Mathematics\\
  University of California, Los Angeles\\
  Los Angeles, CA 90095 \\
  \texttt{xdong99@ucla.edu} \\
   \And
 Koffi Enakoutsa \\
  Department of Mathematics\\
  University of California, Los Angeles\\
  Los Angeles, CA 90095 \\
  \texttt{koffi@math.ucla.edu} \\
}

\begin{document}
\maketitle
\begin{abstract}
    Morrey conjecture deals with two properties of functions which are known as quasi-convexity and rank-one convexity. It is well established that every function satisfying the quasi-convexity property also satisfies rank-one convexity. Morrey (1952) conjectured that the reversed implication will not always hold. In 1992, {\v{S}}ver{\'a}k found a counterexample to prove that Morrey conjecture is true in three dimensional case. The planar case remains, however, open and interesting because of its connections to complex analysis, harmonic analysis, geometric function theory, probability, martingales, differential inclusions and planar non-linear elasticity. Checking analytically these notions is a very difficult task as the quasi-convexity criterion is of non-local type, especially for vector-valued functions. That's why we perform some numerical simulations based on gradient descent algorithms using Dacorogna and Marcellini example functions. Our numerical results indicate that Morrey Conjecture holds true.
\end{abstract}

\keywords{Morrey Conjecture \and Quasi-convexity \and Gradient Descent \and Non-convex Optimization \and Iwaniec Conjecture}

\section{Introduction}
In the 1950's, Charles Morrey worked to find what is the correct notion of convexity in the context of calculus of variations. To address this, he considered the following minimization problem: 
\begin{equation}\label{eqn:minimization problem}
    \text{min}\left\{I(\phi), \phi \in \mathrm{W}^{1,\infty}_0(\Omega, \mathbb{R}^m)\right\} \quad \text{with} \quad I(\phi) = \int_{\Omega} f ( \nabla {\phi} (\mathbf{x})) d \mathbf{x}
\end{equation}
where $\Omega \subset \mathbb{R}^n$ is a bounded open set, $f: \mathbb{R}^{n \times m} \to \mathbb{R}$ is a continuous function, $\mathrm{W}^{1,\infty}_0(\Omega, \mathbb{R}^m)$ denotes the space of $\mathbb{R}^m$-valued Lipschitz functions vanishing on the boundary $\partial\Omega$. Solving the problem in Eq.(\ref{eqn:minimization problem}) is equivalent to proving that the functional $I$ is weakly lower semi-continuous (ellipticity condition) in the Sobolev space $\mathrm{W}^{1,\infty}_0$. The latter is a rather difficult problem and has not yet received a fully satisfactory answer. It was first formulated by Bliss in 1937 in his seminar on the calculus of variations and has received considerable attention, in particular by Albert \cite{albert1938quadratic}, Reid \cite{reid1938theorem}. In addition, MacShane \cite{macshane1939condition}, Hestenes and MacShane \cite{hestenes1940theorem}, Terpstra \cite{terpstra1939}, Van Hove \cite{van1947extension}, Serre \cite{serre1983} and Marcellini \cite{marcellini1984quasiconvex} for the quadratic case and in a more general context by Morrey \cite{morrey1952quasi}, \cite{morrey2009multiple} (see \cite{Ball77}, \cite{ball1987does}, \cite{ball1984w1}, \cite{dacorogna2008} for more details on the quadratic cases). 
$\\$

C. Morrey wanted to define the conditions of convexity on the function $f$ (not including growth or smoothness conditions) that guaranteed that the problem Eq.(\ref{eqn:minimization problem}) is a very well posed variational problem, which translates into the existence of a minimizer for the minimization problem in Eq.(\ref{eqn:minimization problem}). Morrey found out that the functional $I$ in Eq.(\ref{eqn:minimization problem}) is weakly lower semi-continuous if and only if the function $f$ has the quasi-convexity property, which is defined by the Jensen inequality:
\begin{equation}\label{eqn:quasiconvex definition}
    \int_\Omega f(\xi) dx \leq \int_\Omega f(\xi + \nabla \phi(x))dx
\end{equation}
which must be valid for any $\xi \in \mathbb{R}^{n \times m}$ and any $\phi \in \mathrm{W}^{1,\infty}_0(\Omega, \mathbb{R}^m)$.
However, due to its non-local character \cite{kristensen1999} and the presence of vector functions that are gradients in the inequality Eq.(\ref{eqn:quasiconvex definition}),  the quasi-convexity property of a function is a rather difficult criterion to meet. This has motivated Morrey to look for local convexity conditions on $f$ that can guarantee the weakly lower semi-continuity of the functional $I$ in Eq.(\ref{eqn:minimization problem}), see for instance \cite{alibert1992}, \cite{ball1990sets}, \cite{chipot1988equilibrium},  \cite{dacorogna1988counterexample}, \cite{dacorogna1990}, \cite{serre1983formes}, \cite{sivaloganathan1988implications}, \cite{sverak1991quasiconvex}, \cite{vsverak1990examples}, \cite{tartar1979compensated},  and \cite{terpstra1939} as well as the references therein. Some of these conditions are poly-convexity and rank-one convexity. 
$\\$

The notion of poly-convexity was introduced in the context of hyper elasticity theory by Sir John Ball in a pioneering paper \cite{Ball77}. Poly-convexity of a function $f$ is defined as
\begin{equation}\label{eqn:polyconvex definition}
    f(\xi) = \phi(\xi, det\xi)
\end{equation}
for every $\xi \in \mathbb{R}^{2 \times 2}$. Poly-convexity of a strain-energy function is a very important mathematical condition, especially when one tries to solve boundary-value problems.
Also, Alexander Mielke \cite{mielke2005} defined necessary and sufficient conditions for poly-convexity in arbitrary special dimension. Ball \cite{Ball77} suggested that poly-convexity is equivalent to a sufficient condition for quasi-convexity, but there are few known examples of quasi-convex functions that are not poly-convex, establishing that poly-convexity and quasi-convexity conditions are not equivalent. 
$\\$

The other potential local convexity property sought as a replacement for quasi-convexity is rank-one convexity. This condition is defined as 
\begin{equation}\label{eqn:rankoneconvex definition}
    f(\lambda \xi + (1-\lambda)\eta) \leq \lambda f(\xi) + (1-\lambda) f(\eta)
\end{equation}
for every $\lambda \in [0,1]$, $\xi,\eta \in \mathbb{R}^{n \times m}$ with $det(\xi - \eta) = 0$. It is easy to check that the implication
\begin{equation*}\label{eqn:implication}
    f \text{ quasi-convex } \implies f \text{ rank-one convex } 
\end{equation*}
hold for functions $f$ defined on $\mathbb{R}^{n \times m}$. This is done by choosing an appropriate mapping that converts the non-local quasi-convexity condition to the local rank-one convexity condition. The reverse implication does not hold in general for $n > 2$. For instance, a recent remarkable counterexample, in favor of rank-one convexity does not imply quasi-convexity in general, was introduced by Vladimir \v{S}ver\'ak \cite{sverak} and is valid for $m \geq 3$. Whether this latter implication holds for $n = 2$ is still an open question: the conjecture that rank-one convexity and quasi-convexity are not equivalent is also called Morrey Conjecture, see Morrey \cite{morrey1952quasi}. If Morrey Conjecture holds true, many mathematical and/or mechanical modeling methods for material behaviors would have a much more robust theoretical foundations; indeed, for composite materials for instance, the question of whether composites can be constructed with sequential laminates, see \cite{milton2003theory}, would have been resolved as the existence of non quasi-convex but rank-one convex functions is linked to this question. Also, due to their good mathematical structure in terms of variational principles as explained in Gutierrez and Villavicencio \cite{gutierrez2007optimization}, quasi-convex functions are used in modeling phase transition in solids as demonstrated in \cite{carstensen2005ten}, shape optimization (see for instance Pedregal \cite{pedregal2005vector}) , and in fracture mechanics of materials, see Francfort and Marigo \cite{francfort1998revisiting}.   
$\\$

For some classes of functions on $\mathbb{R}^{2 \times 2}$, however, several works have established that the two convexity properties are equivalent, see for instance the works by \cite{muller1999}, \cite{rosakis1994}, \cite{serre1983}, \cite{terpstra1939}. In this spirit, Martin et al. \cite{martin2017} have shown in the context of non-linear elasticity that any energy function $ W: GL^{+} (2) \to \mathbb{R} $ which is isotropic and objective (i.e. bi-$SO(2)$-invariant) as well as isochoric is rank-one convex if and only if the energy density $W$ is poly-convex and doing so gives a negative response to Morrey Conjecture as poly-convexity implies quasi-convexity. In addition, for quadratic types of functions, the equivalence between rank-one and quasi-convexity can be established using Plancherel's formula. The resolution of this equivalence could have big impacts in the theory of conformal mappings for two-component maps case. In particular, if the equivalence between rank-one convexity and quasi-convexity for two component-mappings turns out to be true, then the norm of the corresponding Beurling-Ahlfors transform equals $p^{\star} - 1$, see for instance \cite{iwaniec2002}. It is interesting to point out that Morrey Conjecture also has some connections with the Iwaniec conjecture \cite{iwaniec1982extremal}, the solution of which could impact the resolution of the Morrey Conjecture problem. Indeed, if the equivalence rank-one convexity and quasi-convexity is true, this would imply that the Iwaniec conjecture is true. By Mathematical contraposition, if the Iwaniec conjecture does not hold true, then rank-one convexity would not imply quasi-convexity. Thus, Morrey Conjecture would have been true.  
$\\$

Numerous attempts have been made to construct examples of functions that are rank-one convex, but not quasi-convex \cite{alibert1992}, \cite{dacorogna2008}, \cite{dacorogna1990}. The complexity of the involved calculations has not yet permitted a complete analytical study of such examples, see for instance \cite{duke}, \cite{dacorogna1990}, \cite{gremaud}. Even, the counterexample provided in Vladmir \v{S}ver\'ak \cite{sverak} seems to be a purely three dimension case as many attempts to translate it into a two dimensional setting failed, see \cite{bandeira2011}, \cite{pedregal1996}, \cite{pedregal1998} for references. We could not find other counterexamples in the literature. In a recent note, Pedregal \cite{pedregal2014} provides some evidence in favor of the Morrey Conjecture for two-component maps in dimension two by giving an explicit family of maps parameterized by $\tau$ and proving that for small values of $\tau$ these maps can not be achieved by lamination. As well stated in \cite{pedregal1997}, this will be equivalent to the assertion that there are some rank-one functions that are not quasi-convex, and thereby confirming the validation of Morrey Conjecture. Even though this approach might yield potential counterexample candidates, Pedregal \cite{pedregal2014} recognized himself that the procedure tends to be more involved than in the situation examined by \v{S}ver\'ak. This is a good reason, as a first step, to use numerical analysis to study the problem of whether or not the Morrey Conjecture is valid. One of such works was performed some years ago by Dacorogna et al. \cite{dacorogna1990} on the example of Dacorogna and Marcellini \cite{dacorogna1988counterexample} energy density function; the numerical results of these investigations indicated that Dacorogna and Marcellini example which is rank-one convex, is also quasi-convex.
$\\$

The problem considered in \cite{dacorogna1990} is as follows: 
for $\xi \in \mathbb{R}^{2 \times 2}$ and $\phi \in \mathrm{W}^{1,4}_0(\Omega;\mathbb{R}^2)$, Dacorogna et al. let
\begin{equation}\label{eqg:Eq02}
    f_\gamma(\xi) = \|\xi\|^4 - \gamma \|\xi\|^2 det\xi,
\end{equation}
and
\begin{equation}\label{eqn:Eq03}
    J_\gamma(\xi,\phi) = \int_\Omega[f_\gamma(\xi+\nabla \phi(x)) - f_\gamma(\xi)]d x.
\end{equation}
They choose $\Omega = (0,1) \times (0,1)$ and found that the quasi-convexity of $f_\gamma$ is then equivalent to
\begin{equation}\label{eqn:Eq04}
    \inf_{\xi\in\mathbb{R}^{2 \times 2}}\inf_{\phi \in \mathrm{W}^{1,4}_0(\Omega;\mathbb{R}^2)}\{J_\gamma(\xi,\phi)\} = 0.
\end{equation}
A few remarks are at hand here:
\begin{enumerate}
    \item  First, note that because of the homogeneity of $f_\gamma$, the infimum in Eq.(\ref{eqn:Eq03}) is either $0$ or $-\infty$.
    \item Next, it follows, if $\gamma > \displaystyle\frac{4}{\sqrt{3}}$, then $f_\gamma$ is not rank-one convex and therefore in Eq.(\ref{eqn:Eq03}) the infimum is $-\infty$.
\end{enumerate}

Dacorogna et al. \cite{dacorogna1990} described a numerical approximation to this problem by defining a positive integer $N$ and $h = 1/N$, and a partition $\Omega$ into $\Omega_{ij} = (ih, (i+1)h) \times (jh, (j+1)h), 0 \leq i, j\leq N-1$. Each of these $\Omega_{ij}$ is subdivided into two triangles. They denoted $\tau_k$ this triangulation of $\Omega$ and the triangles by $K$ and they let $P_1, P_2, \ldots, P_M, M = (N-1)^2$, be the internal nodes. Next, they set 
\begin{align*}
    V_h &= \{u\in C^0(\overline{\Omega}): u \text{ is affine on each } K \in \tau_h \text{ and } u=0 \text{ on } \partial\Omega\},\\
    \mathrm{W}_h &= V_h \times V_h \subset \mathrm{W}^{1,4}_0(\Omega;\mathbb{R}^2).
\end{align*}
By fixing $\xi \in \mathbb{R}^{2 \times 2}$ Dacorogna et al. \cite{dacorogna1990} minimize $J_\gamma$ over $\mathrm{W}_h$ using a  gradient descent method, obtained by defining $w^l, l=1,2,\ldots,L$, $d^l=\nabla J_\gamma(w^l)$, $g^l(\alpha) = J_\gamma(w^l+\alpha d^l)$, and updating $w^{l}$ using the explicit gradient update $w^{l+1} = w^l + \overline{\alpha}d^l$, where $\overline{\alpha}$ is obtained by solving $\frac{dg^l}{d\alpha}=0$, using only one step in Newton's method with starting point $\alpha = 0$.
$\\$

The numerical approach used in \cite{dacorogna1990} to solve the above problem is based on a steepest descent algorithm with a crude approximation on the derivation of the gradient of the functional to be minimized. Gremaud \cite{gremaud} used a different numerical approach for the same problem. Unlike in \cite{dacorogna1990}, the corresponding minimizing problem was solved using an annealing-like algorithm. Gremaud's results \cite{gremaud} showed that the example functions considered in \cite{dacorogna1988counterexample} are quasi-convex if and only if they are rank-one convex, contradicting Morrey Conjecture, but confirming Iwaniec conjecture.      
$\\$

Other numerical computation strategies to assess Dacorogna and Marcellini \cite{dacorogna1988counterexample} examples functions with respect to its abilities to provide insights onto the validation/invalidation of Morrey Conjecture exist. Recent works by the authors from Duke University \cite{duke} improved upon the numerical simulations of Dacorogna in an attempt to define a function that is rank-one convex, but not quasi-convex. Duke University's simulations improved on the computational speed and the numerical optimization techniques since the publication of Dacorogna's works. We want to report here also that there has been several numerical attempts to address Morrey Conjecture problem outside of the context of the example proposed by \cite{dacorogna1988counterexample}; among these, let us mention the work of Gutierrez and Villenvicencio \cite{gutierrez2007optimization} where the authors derived an optimization algorithm based on (i) some necessary condition for the quasi-convexity of fourth-degree polynomials that distinguishes between quasi-convex and rank-one convex functions in the three dimensional case, (ii) a characterization of rank-one convex fourth-degree polynomials in terms of infinitely many constraints.
$\\$

The objective of this report is to go somewhat beyond the pioneering works of Dacorogna et al. \cite{dacorogna1990}. We do this by improving on the numerical algorithm these authors used in the gradient descent strategy they proposed. Namely, we calculate the exact expression of the gradient of the functional involved in the optimization problem at hand here and used their approximated values. We solved the minimization problem numerically by using the approximated values of the of the trial functions $\phi$ at each of the nodes of the mesh we used to model the domain $\Omega$. Note here that these values are obtained from an initial trial function $\phi$ that we choose as oscillating functions since in Dacorogna et al.'s numerical computations, these types of functions seem to be promising.  Once the updated values of the trial functions at the nodes are obtained, we used them to check the Jensen's inequality which associated to the quasi-convexity property of the function $f_\gamma$. The initial trial functions to enter the steepest descent iterative algorithm are chosen together with some fixed value of the matrix $\xi$. By randomizing the entries of $\xi$, we successfully used for each of the iteration a new matrix $\xi$. The results indicate that for an appropriate choice of $\gamma$ in the function $f_{\gamma}$ and $\xi$ in Dacorogna \cite{dacorogna1990}, $f_{\gamma}$ is rank-one convex, but the Jensen inequality defining the quasi-convexity is violated, and thereby, for these values of $\gamma$ and $\xi$, $f_{\gamma}$ is rank-one convex, but is not quasi-convex, thus confirming that, at least from the numerical stand point, the Morrey Conjecture holds true. The report is organized as follows:

\begin{itemize}
\item Section \ref{subsec:statement} states the problem to be solved: this problem consists of finding a $\gamma$ in Dacorogna and Marcellini's example \cite{dacorogna1988counterexample} for which $f_\gamma$ is rank-one convex but not quasi-convex. Verifying quasi-convexity of $f_\gamma$ is equivalent to a minimization problem over $\xi \in \mathbb{R}^{2 \times 2}$ and $\phi \in \mathrm{W}^{1,\infty}_0(\Omega, \mathbb{R}^2)$. A non quasi-convexity property on $f_\gamma$ will simply consist of finding some $\xi \in \mathbb{R}^{2 \times 2}$ and $\phi \in \mathrm{W}^{1,\infty}_0(\Omega, \mathbb{R}^2)$ such that Jensen inequality is violated.

\item Section \ref{subsec: numerical approach} presents the numerical approach we use to address the problem considered in the previous section. Our numerical approach is based on a gradient descent algorithm for which, unlike in the previous numerical work on Dacorogna and Marcellini's example \cite{dacorogna1988counterexample}, we calculate the exact gradient of the functional involved in the minimization problem. The functions $\phi$ in the Sobolev space  $\mathrm{W}^{1,\infty}_0$ in our algorithm are approximated by their values at each node of a rectangular grid. These values are updated at each iteration. At the end of an iteration, the values of $\phi$ at each node are used to check Jensen inequality. Also, at each iteration the values of $\xi$ and $\gamma$ are changed to augment our chance to find an "approximated" $\phi$ and a value of $\gamma$ which will violated Jensen inequality.         

\item Finally, Section \ref{sec:results} discusses the results of our numerical simulation. We mainly found that for a given value of $\gamma$ for which $f_{\gamma}$ is rank-one convex we successfully found a function $\phi$ for which Jensen inequality is violated.     
\end{itemize}

\section {Methodology}\label{sec:methodologies}
\subsection{Equations of the Problem}\label{subsec:statement}
In this research, we restrict our attention to Dacorogna and Marcellini's example \cite{dacorogna1988counterexample} to find a $\gamma$ for which $f_\gamma(\xi) = \|\xi\|^4 - \gamma \|\xi\|^2 det\xi$ is rank-one convex but not quasi-convex. Verifying quasi-convexity of $f_\gamma$ consists of checking that: 
\begin{align}\label{eqn:Jensen Inequality}
\begin{split}
    \forall \xi \in \mathbb{R}^{2 \times 2},\forall\phi \in \mathrm{W}^{1,\infty}_0(\Omega, \mathbb{R}^2), J_\gamma(\xi,\phi) = \int_\Omega[f_\gamma(\xi+\nabla \phi(x)) - f_\gamma(\xi)]d x \geq 0
\end{split}
\end{align}

$\\$

Choosing $\Omega = (0,1) \times (0,1)$, the quasi-convexity of $f_\gamma$ is also equivalent to
\begin{equation}\label{eqn:inf}
    \inf_{\xi\in\mathbb{R}^{2 \times 2}}\inf_{\phi \in \mathrm{W}^{1,\infty}_0(\Omega;\mathbb{R}^2)}\{J_\gamma(\xi,\phi)\} = 0.
\end{equation}
$\\$

Non quasi-convexity property of $f_\gamma$ will occur whenever we can find a $\xi \in \mathbb{R}^{2 \times 2}$ and a $\phi \in \mathrm{W}^{1,\infty}_0(\Omega, \mathbb{R}^2)$ such that inequality in Eq.(\ref{eqn:Jensen Inequality}) is violated. A reformulation of this problem was proposed by Cheng et.al. \cite{duke} at Duke University. It consists of looking for a $\phi$ such that
\begin{equation}
    \gamma^* = \sup_{\phi\in \mathrm{W}_0^{1,\infty}(\Omega, \mathbb{R})}\left\{ \frac{\int_\Omega \|\xi+\nabla\phi\|^4-\|\xi\|^4d\Omega}{\int_\Omega\|\xi+\nabla\phi\|^2det(\xi+\nabla\phi)-\|\xi\|^2det\xi d\Omega}\right\}.
\end{equation}
$\\$

If the $\gamma^*$ is found such that $f_\gamma$ is rank-one convex, then Cheng et.al. would have found an example of a function that is rank-one convex but not quasi-convex. In Cheng et.al. approach, $\phi$ is given by a sinusoidal family of vector functions. We didn't follow this approach in this work. However, we recently implemented Cheng et.al. approach and obtained a significant improvement of the results by Cheng et.al., see Dong and Enakoutsa \cite{dong2022some}.

\subsection{Numerical Implementation}\label{subsec: numerical approach}
\subsubsection{The Algorithm}\label{subsubsec: Algorithm}

The objective of our numerical approach is to find some suitable $\phi$, $\xi$, and $\gamma$ to construct a counterexample function based on Dacorogna and Marcellini's example function \cite{dacorogna1988counterexample} which is rank-one convex but not quasi-convex. The chances of finding these three quantities will be maximized if we can rich a generate and cover different permutations of $\phi$, $\xi$, and $\gamma$. To do this, we perform a steepest descent algorithm on the function $\phi$ for some value of the matrix $\xi$ using MATLAB$\copyright$ software \cite{MATLAB:2021}. 
$\\$

\begin{figure}[htp]
    \centering
    \includegraphics[width=8cm]{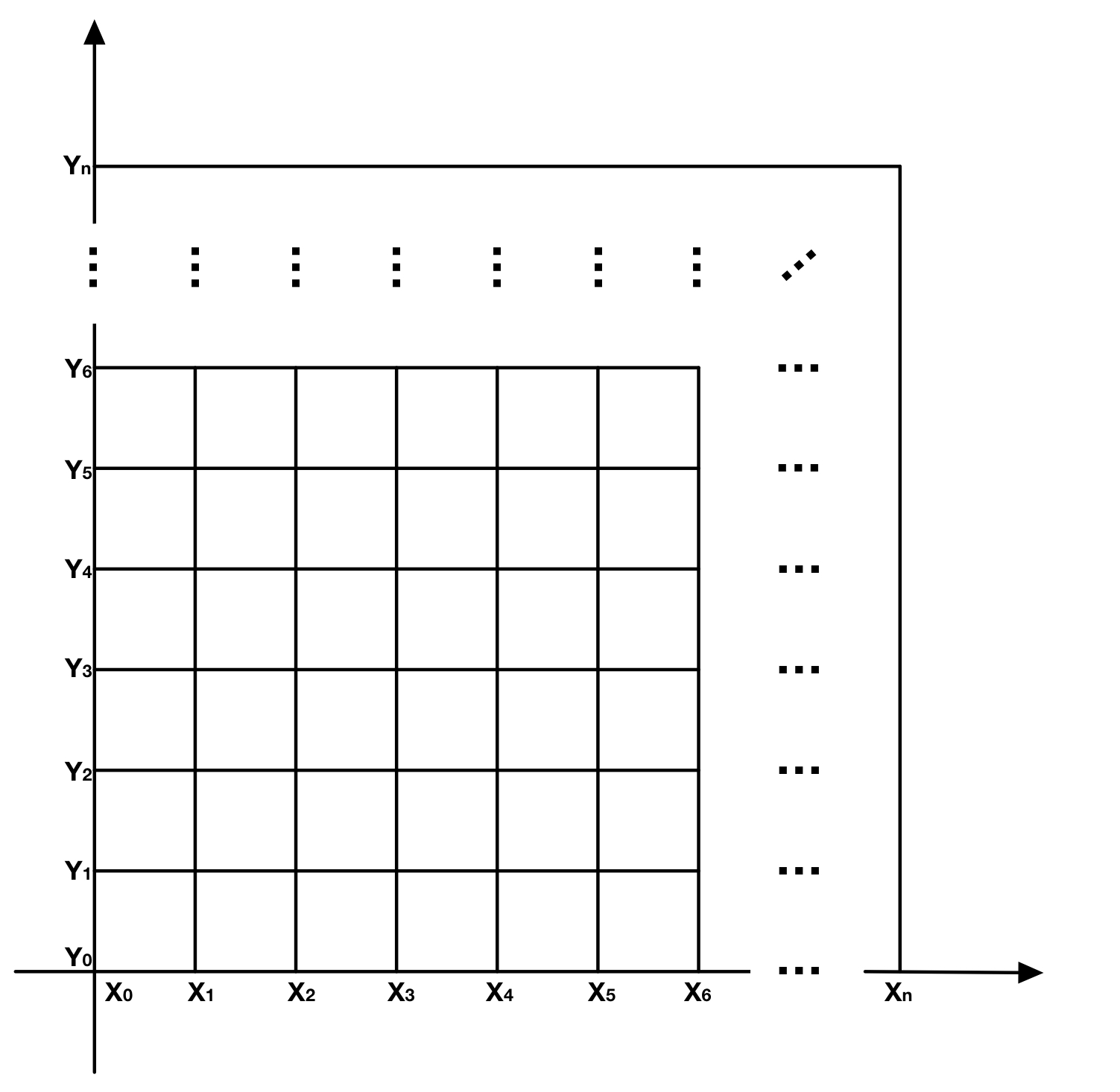}
    \caption{Mesh of the Square Domain $\Omega$}
    \label{MeshGraph}
\end{figure}

First, we create a rectangular grid to model $\Omega = (0,1) \times (0,1)$ and store the values of the coordinates of the nodes in the grid. The mesh size of each square element in the grid is $\displaystyle\frac{1}{n}$ where $n$ is the number of subintervals on each of the axes, see Figure \ref{MeshGraph}.
$\\$

Next, we fix $\xi \in \mathbb{R}^{2\times 2}$ and check the Jensen inequality $J_\gamma(\xi, \phi)$ for the approximated values for function $\phi$ at all the nodes. The values of $J_\gamma(\xi, \phi)$ are computed by a numerical double integral over $\Omega$ of the integrand $f_\gamma(\xi+\nabla\phi^k)-f_\gamma(\xi)$ where $\phi^k$ is updated using the line search recursive formula
\begin{equation}\label{eqn:recursive formula}
    \phi^{k+1} = \phi^k - \tau^k g^k
\end{equation} at all the nodes. 
$\\$

In Eq.(\ref{eqn:recursive formula}), $g^k = \delta J_\gamma(\xi,\phi^k)$ is obtained using G\^ateaux differentiation of $J_\gamma(\xi,\phi^k)$ with respect to $\phi^k$ as
    \begin{equation}\label{eqn:delta J}
        \delta J_\gamma(\xi,\phi^k) = - \displaystyle\sum_{i,j=1}^2\displaystyle\frac{\partial}{\partial x_j}\displaystyle\frac{\partial}{\partial \xi_{ij}}f_\gamma(\xi + \nabla\phi^k).
    \end{equation} Details of this derivation of the formula on the right side of Eq.(\ref{eqn:delta J}) can be found in the \hyperref[app:calculation grad]{Appendix}. With Eq.(\ref{eqn:delta J}), Eq.(\ref{eqn:recursive formula}) becomes
    \begin{center}
    \begin{equation}
    \boxed{\begin{aligned}\label{update algorithm matlab} 
    \vspace{0.25cm}
    \phi_1^{k+1} &= \phi_1^k + \tau^k\displaystyle\sum_{j=1}^2\displaystyle\frac{\partial}{\partial x_j}\displaystyle\frac{\partial}{\partial \xi_{1j}}f_\gamma(\xi + \nabla\phi^k)\\
    \vspace{0.25cm}
    &= \phi_1^k + \tau^k\displaystyle(\displaystyle\frac{\partial}{\partial x_1}\displaystyle\frac{\partial}{\partial \xi_{11}}f_\gamma(\xi + \nabla\phi^k) + \frac{\partial}{\partial x_2}\frac{\partial}{\partial \xi_{12}}f_\gamma(\xi + \nabla\phi^k))\\
    \vspace{0.25cm}
    \phi_2^{k+1} &= \phi_2^k + \tau^k\displaystyle\sum_{j=1}^2\displaystyle\frac{\partial}{\partial x_j}\displaystyle\frac{\partial}{\partial \xi_{1j}}f_\gamma(\xi + \nabla\phi^k)\\
    &= \phi_2^k + \tau^k(\displaystyle\frac{\partial}{\partial x_1}\displaystyle\frac{\partial}{\partial \xi_{12}}f_\gamma(\xi + \nabla\phi^k) + \displaystyle\frac{\partial}{\partial x_2}\displaystyle\frac{\partial}{\partial \xi_{22}}f_\gamma(\xi + \nabla\phi^k))
\end{aligned}}
\end{equation}
\end{center}
$\\$
after expansion.
$\\$

The second term in the recursive formula Eq.(\ref{eqn:recursive formula}), $\tau^k$ is defined as in the steepest gradient descent algorithm formulation:
    \begin{equation}\label{eqn:tau}
        \tau^k = \underset{\alpha > 0}{\mathrm{argmin}}\ J_\gamma(\phi^k - \alpha g^k).
    \end{equation}
    
Let's introduce a new function $\psi(\alpha)$ such that
\begin{align}
\begin{split}
    \psi(\alpha) \coloneqq J_\gamma(\phi^k - \alpha g^k).
\end{split}
\end{align}

Then finding $\tau^k$ in Eq.(\ref{eqn:tau}) is equivalent to finding $\alpha$ such that
\begin{align}
\begin{split}
    \frac{d}{d\alpha}\psi(\alpha) = 0\\
\end{split}
\end{align}

Expanding $\psi(\alpha)$, we get
\begin{align}\label{eqn:expanding dalpha}
\begin{split}
    \frac{d}{d\alpha}\int_\Omega f_\gamma(\xi+\nabla(\phi^k - \alpha\delta J_\gamma(\phi^k)))d\Omega = 0
\end{split}
\end{align}

Taking the derivative in Eq.(\ref{eqn:expanding dalpha}) we get
\begin{align}\label{eqn:inner product}
\begin{split}
         \displaystyle\frac{d}{d\alpha}\psi(\alpha)
         = \int_\Omega \left\langle Df_\gamma(\xi+\nabla(\phi^k - \alpha\delta J_\gamma(\phi^k)), \frac{\partial}{\partial\alpha}(\nabla(\phi^k - \alpha \delta J_\gamma(\phi^k)))\right\rangle d\Omega
\end{split}
\end{align}
where $\langle\ ,\ \rangle$ is the inner product defined on $\mathrm{W}^{1,4}_0(\Omega;\mathbb{R}^2)$. As a consequence, Eq.(\ref{eqn:inner product}) results in
\begin{align}\label{eqn:after expanding}
\begin{split}
    \displaystyle\frac{d}{d\alpha}\psi(\alpha) = \int_\Omega \left\langle \delta J_\gamma(\phi^k - \alpha \delta J_\gamma(\phi^k)), \delta J_\gamma(\phi^k)\right\rangle d\Omega.
\end{split}
\end{align}

Assume that $h(\alpha) = \displaystyle\frac{d}{d\alpha}\psi(\alpha)$; due to the difficulties that arose to obtain an analytic expression for the root of the equation $h(\alpha) = 0$, the structure of Eq.(\ref{eqn:after expanding}) enables us to use the secant method instead of other root finding methods such as the Newton or bisection methods. The iterative sequence of the secant method is given by the standard formula
\begin{align}
\begin{split}
    \alpha^{n+1} = \alpha^n - h(\alpha^n)\displaystyle\frac{\alpha^n-\alpha^{n-1}}{h(\alpha^n)-h(\alpha^{n-1})},
\end{split}
\end{align}
with $\alpha^0, \alpha^1$ as two initial guess real numbers.
$\\$

A flow chart of the algorithm is described as below:
$\\$

\begin{center}
\noindent\fbox{%
\begin{varwidth}{\dimexpr\linewidth-2\fboxsep-2\fboxrule\relax}
\begin{algorithmic}[1]\label{algorithm}
\State Create a rectangular grid and store the coordinates of the nodes in the grid
\State Choose $2< \gamma < \displaystyle\frac{4}{\sqrt{3}} \text{ and } \xi \in \mathbb{R}^{2 \times 2}$ 
\State Initialize $k=0$ and $\phi^0$ by finding its value at all nodes
\State Calculate $g^k = \delta J_\gamma(\phi^k)$ at all the nodes
\State Update the values of function $\phi$ using the recursive formula: \\
$\phi^{k+1} = \phi^k - \tau^k g^k$ at all the nodes
\State Set the values of $\phi$ on boundary nodes to $0$
\State Compute the step size $\tau^k = \underset{\alpha > 0}{\mathrm{argmin}}\ J_\gamma(\phi^k - \alpha g^k)$ using the secant method.
\If {$J_\gamma(\xi, \phi^{k+1}) > 0$\ (quasi-convexity)}
    \State Set $\gamma \Leftarrow \gamma - \delta\gamma$, $\phi^k \Leftarrow \phi^{k+1}$, $k \Leftarrow k+1$
    \State Go to step 4 (till $\gamma = 2$)
\Else
    \State Repeat Step 10 and Step 11
\EndIf
\end{algorithmic}
\end{varwidth}%
}
\end{center}

\section{Numerical Results}\label{sec:results}
This section presents the numerical results we obtained using the algorithm described in the previous section. We shall distinguish two cases. In the first case we assume that $\xi$ is fixed. In the second type of simulations, the matrix $\xi$ is changed from one iteration to another one thanks to the MATLAB $\copyright$ \cite{MATLAB:2021} function $RAND$, a random number generator function, which changes the values of the entries in the matrix $\xi$ at each iteration.
\subsection{Numerical Results when $\xi$ is Fixed during the Gradient Descent Iterations}
In this case, we fix $\xi$ to $\xi = [1,0; 0,\sqrt{3}]$ as in Dacorogna et.al. \cite{dacorogna1990}, and change the values of $\gamma$ starting from $\displaystyle\frac{4}{\sqrt{3}}$ to $2$. For each value of $\gamma$, we update $\phi^k$ using the recursive formulas in Eq.(\ref{update algorithm matlab}) and Eq.(\ref{eqn:tau}) of the steepest descent algorithm. Then we check Jensen inequality $J_\gamma(\xi, \phi)$ by computing the double integral over $\Omega$ of the integrand $f_\gamma(\xi+\nabla\phi^k)-f_\gamma(\xi)$. 
$\\$

\begin{figure}[htp]\label{100sqrt3}
  \centering
  \begin{subfigure}[t]{.45\textwidth}
    \centering
    \includegraphics[width=\linewidth]{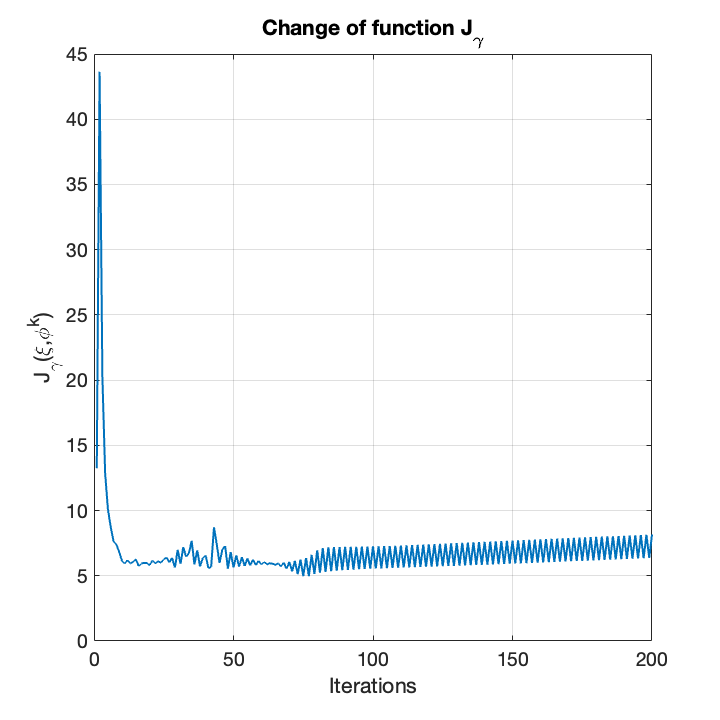}
    \caption{Change of function during the steepest descent on $\phi_1^0(x,y) = [
    \displaystyle\frac{1}{2\pi}\sin(2\pi x), \displaystyle\frac{1}{2\pi}\sin(2\pi y)
 ]$ with fixed $\xi = [0,0; 0,0]$ and changing $\gamma$}
  \end{subfigure}
  \hfill
  \begin{subfigure}[t]{.45\textwidth}
    \centering
    \includegraphics[width=\linewidth]{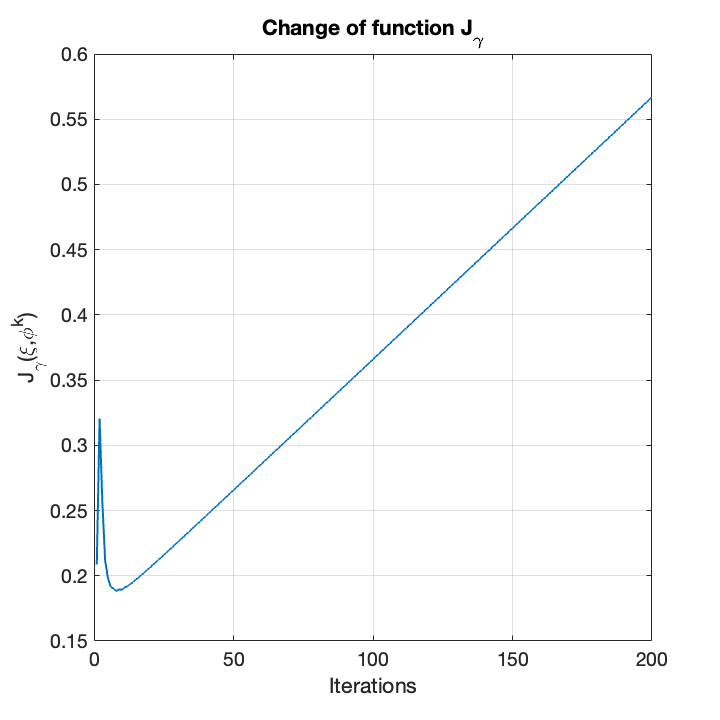}
    \caption{Change of function during the steepest descent on $\phi_2^0(x,y) = [
    \sin(x(x-1)y(y-1)), \sin(x(x-1)y(y-1))^2
    ]$ with fixed $\xi = [0,0; 0,0]$ and changing $\gamma$}
  \end{subfigure}
  \caption{Change of function during the steepest descent with fixed $\xi = [1,0; 0,\sqrt{3}]$ and changing $\gamma$}
\end{figure}

\hyperref[100sqrt3]{Figure 2(a)} represents the change of values of $J_\gamma(\xi, \phi)$ per iteration when the initial guess function to begin the steepest descent algorithm is: $\phi_1^0(x,y) = \begin{bmatrix}
    \displaystyle\frac{1}{2\pi}\sin(2\pi x), \displaystyle\frac{1}{2\pi}\sin(2\pi y)
\end{bmatrix}$. The figure shows that the change of $J_\gamma(\xi, \phi)$ is falling until it reaches a value where it almost becomes constant (the oscillations around such plateau value can be explained by the numerical approximation we used to calculate the gradient $\displaystyle\frac{\partial}{\partial x_j}\frac{\partial}{\partial \xi_{ij}}$ over the square grid). Let's note that the value of $J_\gamma(\xi, \phi)$ is descending but it remains positive for all of the $\phi^k$ and $\gamma$ we used. Therefore, the Jensen inequality cannot be violated (the function remains quasi-convex while for those values of $\gamma$ it is rank-one convex).
$\\$

\hyperref[100sqrt3]{Figure 2(b)} is the analogous of \hyperref[FixedXiP1]{Figure 2} for $\phi_2^0(x,y) = \begin{bmatrix}
    \sin(x(x-1)y(y-1)), \sin(x(x-1)y(y-1))^2
\end{bmatrix}$. Here, there is a quick descent of the function $J_\gamma(\xi, \phi)$ which starts to rise after a minimal value. Note that, like in the first simulation with $\phi_1^0$, $J_\gamma(\xi, \phi)$ remains positive for all of the $\phi^k$ and $\gamma$ we used. 

Next, we fix $\xi$ to $[0,0; 0,0]$ as in Dacorogna et.al. \cite{dacorogna1990} and Cheng et.al. \cite{duke}. For this case, we used four different initial guess functions $\phi^0$ for the simulations:
\begin{equation}
  \phi^0 (x,y) =
    \begin{cases}
      & \phi^0_1(x,y) = \begin{bmatrix}
    \sin(x(x-1)y(y-1)), \sin(x(x-1)y(y-1))^2
    \end{bmatrix}\\\
      & \phi^0_2(x,y) = \begin{bmatrix}
    (x(x-1)y(y-1)), 0
    \end{bmatrix} \\
      & \phi^0_3(x,y) = \begin{bmatrix}
    \displaystyle\frac{1}{2\pi}\sin(2\pi x), \displaystyle\frac{1}{2\pi}\sin(2\pi y)
    \end{bmatrix} \\
    & \phi^0_4(x,y) = \begin{bmatrix}
    \displaystyle\frac{1}{100}\sin(\pi x), \displaystyle\frac{1}{100}\sin(\displaystyle\frac{3\pi}{2}y)
    \end{bmatrix}
    \end{cases}       
\end{equation}

\begin{figure}[htp]
  \centering
  \begin{subfigure}[t]{.45\textwidth}
    \centering
    \includegraphics[width=\linewidth]{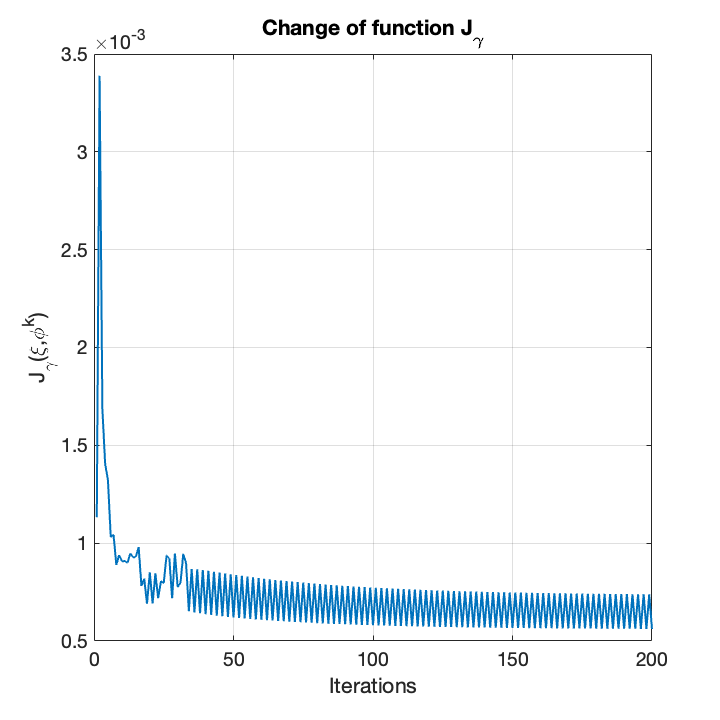}
    \caption{Change of function during the steepest descent on $\phi_1(x,y) = [
    \sin(x(x-1)y(y-1)), \sin(x(x-1)y(y-1))^2
    ]$ with fixed $\xi = [0,0; 0,0]$ and changing $\gamma$}
  \end{subfigure}
  \hfill
  \begin{subfigure}[t]{.45\textwidth}
    \centering
    \includegraphics[width=\linewidth]{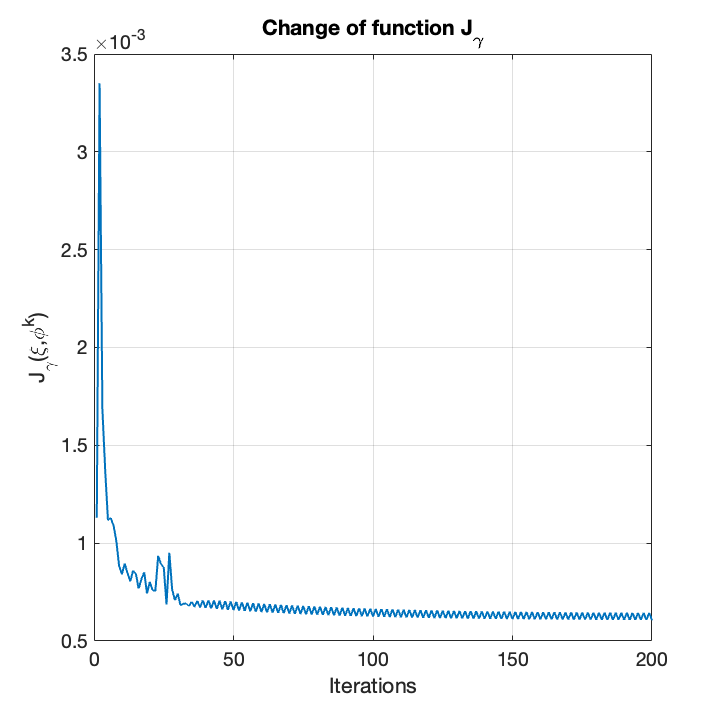}
    \caption{Change of function during the steepest descent on $\phi^0_2(x,y) = [
    (x(x-1)y(y-1)), 0
    ]$ with fixed $\xi = [0,0; 0,0]$ and changing $\gamma$}
  \end{subfigure}

  \medskip

  \begin{subfigure}[t]{.45\textwidth}
    \centering
    \includegraphics[width=\linewidth]{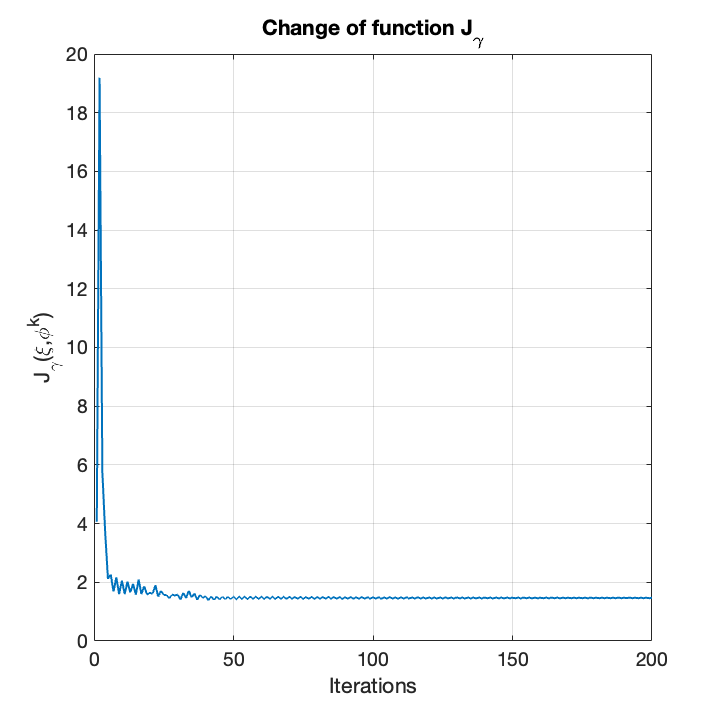}
    \caption{Change of function during the steepest descent on $\phi^0_3(x,y) = [
    \displaystyle\frac{1}{2\pi}\sin(2\pi x), \displaystyle\frac{1}{2\pi}\sin(2\pi y)
    ]$ with fixed $\xi = [0,0; 0,0]$ and changing $\gamma$}
  \end{subfigure}
  \hfill
  \begin{subfigure}[t]{.45\textwidth}
    \centering
    \includegraphics[width=\linewidth]{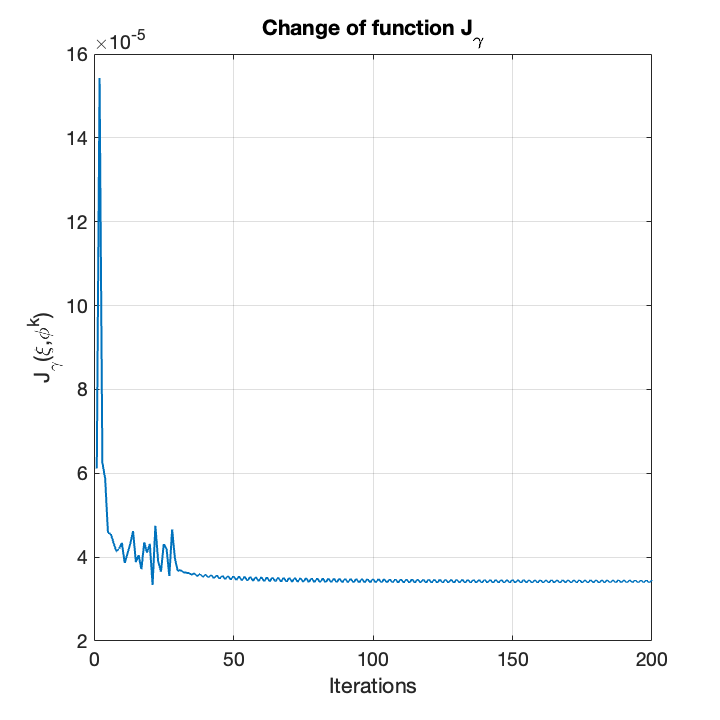}
    \caption{Change of function during the steepest descent on $\phi^0_4(x,y) = [
    \displaystyle\frac{1}{100}\sin(\pi x), \displaystyle\frac{1}{100}\sin(\displaystyle\frac{3\pi}{2}y)
    ]$ with fixed $\xi = [0,0; 0,0]$ and changing $\gamma$}
  \end{subfigure}
  \caption{Change of function during the steepest descent with fixed $\xi = [0,0; 0,0]$ and changing $\gamma$}
  \label{0000 figures}
\end{figure}

The figures (\hyperref[0000 figures]{Figure 3(a) - (d)}) obtained for each of the initial guess has overall same behaviors: the curves fall sharply until a minimal value where it does not change as the iteration increases. Here also, the oscillations obtained after the minimal values may be due to the numerical approximations we made to find the derivatives involved in the gradients to update the steepest descent iteration. For the four example initial functions we studied, the figures (\hyperref[0000 figures]{Figure 3(a) - (d)}) demonstrate that Jensen inequality was not violated (the functions $f_\gamma$ remain quasi-convex). These simulations seem to illustrate the findings of some of the work \cite{dacorogna1990}, \cite{voss2022numerical}, \cite{gremaud} who state that rank-one convexity implies quasi-convexity. While we know such implication will yield interesting results such as finding sharp constants for the $l^p$ bound of the bound of Beurling Alfhors transform (or the complex Hilbert transform) \cite{iwaniec1982extremal}, we must say that because the number of these initial guesses is small, we can not conclude that the methodology we are using support the validity of the previous implication. This is a good reason to seek for other approaches with our algorithm to find a suitable example that can illustrate Morrey Conjecture. In this way, we decide, while keeping the algorithm the same, to vary the values of the entries of $\xi$ at each of the iterations of our steepest gradient descent method. Our next section will report the results of such strategy. 

\subsection{Numerical Results when $\xi$ is changing during the Gradient Descent Iterations}
In the results that follow, we modify the \hyperref[algorithm]{algorithm} by replacing the fixed $\xi$ at each of the iterations at the steepest gradient descent method by changing its values at each of the iterations using the built-in $RAND$ function of MATLAB$\copyright$ \cite{MATLAB:2021}. Doing so allows us to enlarge the set of $\xi$ that can be tested so that Jensen inequality is violated. Our numerical simulations with the modified algorithm reveal that for three of the four initial guess functions we used previously, there are sets of $\phi^k$, $\gamma$, and $\xi$ for which Dacorogna's $f_\gamma$ is rank-one convex but the Jensen inequality is violated, which translates into Dacorogna \cite{dacorogna1988counterexample} example function is rank-one convex but not quasi-convex. For instance, see \hyperref[Table 1]{Table 1}, for 
\begin{equation}
\phi_1^0 (x, y) = \begin{bmatrix}
\sin(x(x-1)y(y-1)), \quad \sin(x(x-1)y(y-1))^2
\end{bmatrix}
\end{equation}
our simulations indicate that for $\xi = \begin{bmatrix}
0.79982, 0.099004;
0.20266, 0.83511
\end{bmatrix}$ and $\gamma = 2.2723$, the Jensen integral is giving a negative value: $J_\gamma(\xi,\phi^k) = -0.018949$. Also, for the case where 
\begin{equation}
\phi_2^0 (x, y) = 
\begin{bmatrix}
    (x(x-1)y(y-1)),\quad 0
\end{bmatrix}
\end{equation}
our simulations show that for $\xi = \begin{bmatrix}
0.95184, 0.28623;
0.083266, 0.75639
\end{bmatrix}$
and $\gamma = 2.298$, the Jensen integral is giving a negative value: $J_\gamma(\xi,\phi^k) = -0.009052$ as shown in \hyperref[Table 2]{Table 2}.
$\\$

Analogous results are obtained for 
\begin{equation}
\phi_4^0 (x, y) = \begin{bmatrix}
    \displaystyle\frac{1}{100}\sin(\pi x), \quad \displaystyle\frac{1}{100}\sin(\frac{3\pi}{2}y)
\end{bmatrix}
\end{equation}
see \hyperref[fig:randEpsFun4]{Figure 7} and \hyperref[Table 3]{Table 3}. The negative values obtained for the Jensen function seem to be small in the three cases, this might be due to the magnitude of the entries of $\xi$ as the $RAND$ function we are using are for entries that are positive and less than 1, so we repeat another round of simulations with the same modified algorithm, but this time with higher order of magnitude of the entries of $\xi$. Same behaviors as previously were observed: we found suitable $\phi^k$, $\gamma$, and $\xi$ for which Morrey Conjecture holds, see \hyperref[fig:scaledXi]{Figure 8} and \hyperref[Table 4]{Table 4}.
$\\$

Similarly, all the three results were obtained with a mesh size $h=1/10$. For the sake of verification, we choose a mesh size of $h = 1/30$, smaller than the size in the previous set of simulations. \hyperref[fig:smallermesh]{Figure 9} and \hyperref[Table 5]{Table 5} present analogous results as for the $1/10$ mesh size: for $\gamma = 2.2679 $ and $\xi = \begin{bmatrix}
    0.97238 & 0.22812\\
    0.044135&0.79754
\end{bmatrix}$, the numerical simulations result reveals that $J_\gamma(\xi,\phi) = -0.006274$.

\subsubsection{$h = \frac{1}{10}$, $\xi$ with all its entries normal distributed between $0$ and $1$}
\begin{enumerate}
    \item For \begin{equation*}
        \phi_1^0(x,y) = \begin{bmatrix}
    \sin(x(x-1)y(y-1)), \sin(x(x-1)y(y-1))^2
\end{bmatrix}
    \end{equation*}

\begin{figure}[htp]\label{fig:randEpsFun1}
    \centering
    \includegraphics[width=12cm]{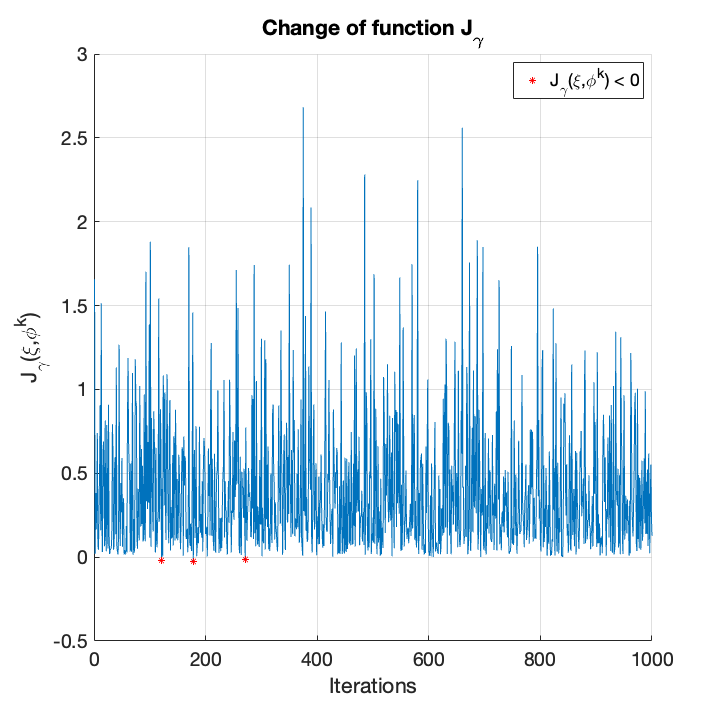}
    \caption{Change of function during the steepest descent on $\phi_1(x,y) = [
    \sin(x(x-1)y(y-1)), \sin(x(x-1)y(y-1))^2
    ]$ with changing $\xi$ and $\gamma$}
\end{figure}
\clearpage
\begin{table}\label{Table 1}
\centering
\begin{tabular}{lllll}
\toprule
    & $\qquad\qquad\xi$ \qquad\qquad\qquad\qquad\qquad & $\quad\gamma$ \qquad\qquad\qquad\qquad& $J_\gamma(\xi, \phi_k)$  \\
\midrule
\\
1. & $\begin{bmatrix}
0.79982 & 0.099004\\
0.20266 & 0.83511
\end{bmatrix}$ & 2.2723 & -0.018949
\\
\\
\midrule
\\
2. &$\begin{bmatrix} 
0.99792 & 0.13205 \\ 
0.028506 & 0.75162
\end{bmatrix}$ & 2.2546  & -0.028785
\\
\\
\midrule
\\
3. &$\begin{bmatrix}  
0.89121 & 0.092286 \\ 
0.05108 & 0.68636
\end{bmatrix}$ & 2.2259 & -0.012748
\\
\\
\bottomrule
\end{tabular}
\caption{Numerical values of $\xi, \gamma, \text{ and } J_\gamma(\xi, \phi_k) \text{ when } J_\gamma(\xi, \phi_k) < 0$}
\end{table}

\clearpage

\item For \begin{equation*}
    \phi_2^0(x,y) = \begin{bmatrix}
    (x(x-1)y(y-1)), 0
\end{bmatrix}
\end{equation*}

\begin{figure}[htp]\label{fig:randXiFun2}
    \centering
    \includegraphics[width=12cm]{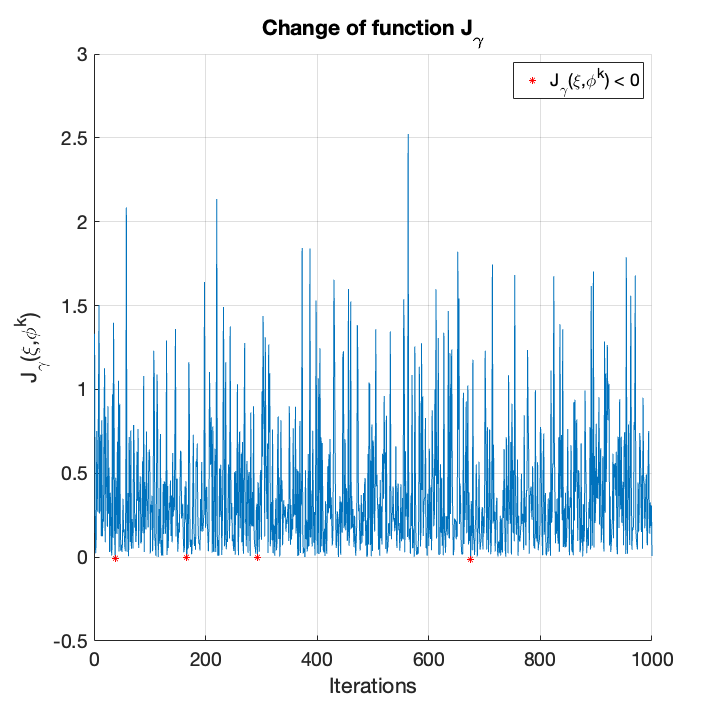}
    \caption{Change of function during the steepest descent on $\phi^0_2(x,y) = [
    (x(x-1)y(y-1)), 0
    ]$ with changing $\xi$ and $\gamma$}
\end{figure}
\clearpage

\begin{table}\label{Table 2}
\centering
\begin{tabular}{llll}
\toprule
    & $\qquad\qquad\xi$ \qquad\qquad\qquad\qquad\qquad & $\quad\gamma$ \qquad\qquad\qquad\qquad& $J_\gamma(\xi, \phi_k)$ \\
\midrule
\\
1. & $\begin{bmatrix}
0.95184 & 0.28623\\
0.083266 & 0.75639
\end{bmatrix}$ & 2.298 & -0.009052
\\
\\
\midrule
\\
2. &$\begin{bmatrix} 
0.57963 & 0.035939 \\ 
0.13491 & 0.64033
\end{bmatrix}$ & 2.2583  & -0.0024306
\\
\\
\midrule
\\
3. &$\begin{bmatrix}  
0.92504 & 0.10624 \\ 
0.18151 & 0.69145
\end{bmatrix}$ & 2.2191 & -0.0025236     
\\
\\
\midrule
\\
4. & $\begin{bmatrix}  
0.9393 & 0.030814 \\ 
0.093527 & 0.92682
\end{bmatrix}$ & 2.1009 & -0.015936   
\\
\\
\bottomrule
\end{tabular}
\caption{Numerical values of $\xi, \gamma, \text{ and } J_\gamma(\xi, \phi_k) \text{ when } J_\gamma(\xi, \phi_k) < 0$}
\end{table}
\clearpage

\item For \begin{equation*}
    \phi_3^0(x,y) = \begin{bmatrix}
    \displaystyle\frac{1}{100}\sin(\pi x), \displaystyle\frac{1}{100}\sin(\frac{3\pi}{2}y)
\end{bmatrix}
\end{equation*}

\begin{figure}[htp]\label{fig:randEpsFun4}
    \centering
    \includegraphics[width=12cm]{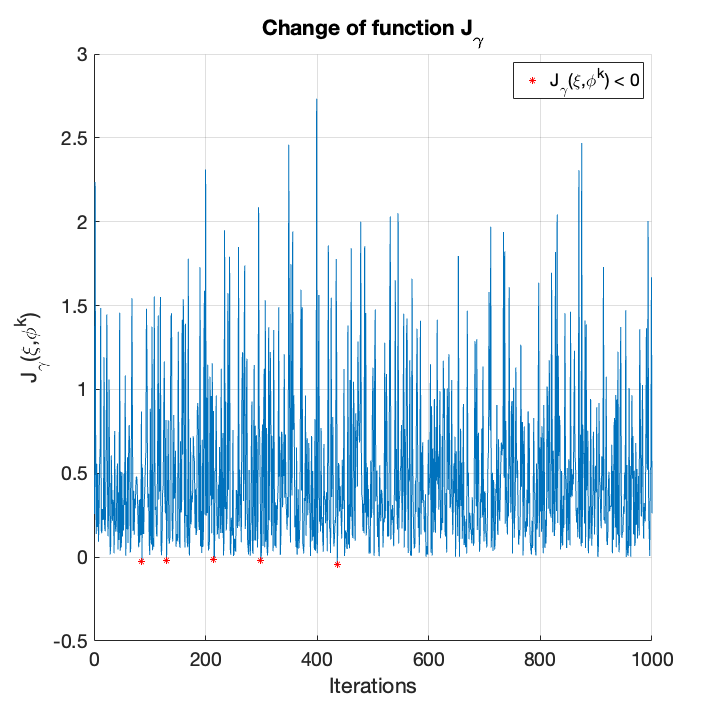}
    \caption{Change of function during the steepest descent on $\phi^0_3(x,y) = [
    \displaystyle\frac{1}{100}\sin(\pi x), \displaystyle\frac{1}{100}\sin(\displaystyle\frac{3\pi}{2}y)
    ]$ with changing $\xi$ and $\gamma$}
\end{figure}
\newpage

\begin{table}\label{Table 3}
\centering
\begin{tabular}{llll}
\toprule
    & $\qquad\qquad\xi$ \qquad\qquad\qquad\qquad\qquad & $\quad\gamma$ \qquad\qquad\qquad\qquad& $J_\gamma(\xi, \phi_k)$  \\
\midrule
\\
1. & $\begin{bmatrix}
0.81147 & 0.26021\\
0.011742 & 0.96949
\end{bmatrix}$ & 2.2831 & -0.027452
\\
\\
\midrule
\\
2. &$\begin{bmatrix} 
0.95366 & 0.26272 \\ 
0.058398 & 0.90193
\end{bmatrix}$ & 2.2695  & -0.022328
\\
\\
\midrule
\\
3. &$\begin{bmatrix}  
0.73263 & 0.11277 \\ 
0.016989 & 0.66646
\end{bmatrix}$ & 2.2432 & -0.011148
\\
\\
4. &$\begin{bmatrix} 
0.81183 & 0.018199 \\ 
0.19432 & 0.94954
\end{bmatrix}$ & 2.2172  & -0.021571
\\
\\
\midrule
\\
5. &$\begin{bmatrix}  
0.96109 & 0.035623 \\ 
0.029006 & 0.94156
\end{bmatrix}$ & 2.1748 & -0.04477
\\
\\
\bottomrule
\end{tabular}
\caption{Numerical values of $\xi, \gamma, \text{ and } J_\gamma(\xi, \phi_k) \text{ when } J_\gamma(\xi, \phi_k) < 0$}
\end{table}
\end{enumerate}
\clearpage
\subsubsection{$h = \frac{1}{20}$, $\xi$ with all its entries normal distributed between $0$ and $1$}
\begin{enumerate}
    \item For \begin{equation*}
        \phi_1^0(x,y) = \begin{bmatrix}
    \sin(x(x-1)y(y-1)), \sin(x(x-1)y(y-1))^2
\end{bmatrix}
    \end{equation*}
\begin{figure}[htp]\label{fig:smallermeshXiFun1}
    \centering
    \includegraphics[width=12cm]{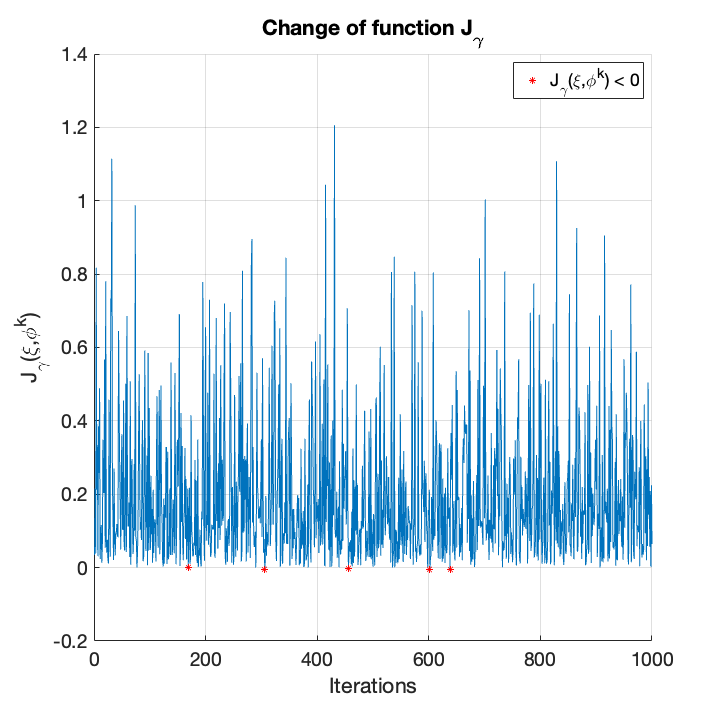}
    \caption{Change of function during the steepest descent on $\phi^0_1(x,y) = [
    \sin(x(x-1)y(y-1)), \sin(x(x-1)y(y-1))^2
    ]$ with changing $\xi$ and $\gamma$}
\end{figure}
\newpage

\begin{table}\label{Table 4}
\centering
\begin{tabular}{llll}
\toprule
    & $\qquad\qquad\xi$ \qquad\qquad\qquad\qquad\qquad & $\quad\gamma$ \qquad\qquad\qquad\qquad& $J_\gamma(\xi, \phi_k)$  \\
\midrule
\\
1. & $\begin{bmatrix}
0.9462 & 0.069215\\
0.13961 & 0.5878
\end{bmatrix}$ & 2.2574  & -0.00032465
\\
\\
\midrule
\\
2. &$\begin{bmatrix} 
0.94937 & 0.041621 \\ 
0.2072 & 0.70936
\end{bmatrix}$ & 2.215  & -0.0062741
\\
\\
\midrule
\\
3. &$\begin{bmatrix}  
0.78787 & 0.046076 \\ 
0.023873 & 0.97483
\end{bmatrix}$ & 2.1686 & -0.0019017
\\
\\
\midrule
\\
4. &$\begin{bmatrix}  
0.98475 & 0.0044322 \\ 
0.072928 & 0.79571
\end{bmatrix}$ & 2.1235 & -0.0054715
\\
\\
\midrule
\\
5. &$\begin{bmatrix}  
0.94015 & 0.08685 \\ 
0.11937 & 0.82277
\end{bmatrix}$ & 2.1123 & -0.0039773
\\
\\
\bottomrule
\end{tabular}
\caption{Numerical values of $\xi, \gamma, \text{ and } J_\gamma(\xi, \phi_k) \text{ when } J_\gamma(\xi, \phi_k) < 0$}
\end{table}
\clearpage

\item For \begin{equation*}
    \phi_2^0(x,y) = \begin{bmatrix}
    (x(x-1)y(y-1)), 0
\end{bmatrix}
\end{equation*}

\begin{figure}[htp]\label{fig:smallermeshXiFun2}
    \centering
    \includegraphics[width=12cm]{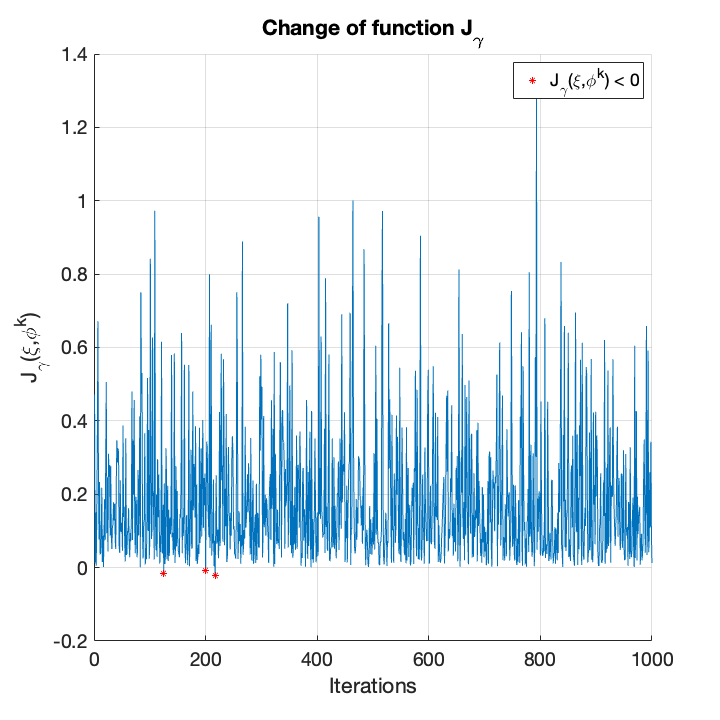}
    \caption{Change of function during the steepest descent on $\phi^0_2(x,y) = [
    (x(x-1)y(y-1)), 0
    ]$ with changing $\xi$ and $\gamma$}
\end{figure}
\newpage

\begin{table}\label{Table 5}
\centering
\begin{tabular}{llll}
\toprule
    & $\qquad\qquad\xi$ \qquad\qquad\qquad\qquad\qquad & $\quad\gamma$ \qquad\qquad\qquad\qquad& $J_\gamma(\xi, \phi_k)$  \\
\midrule
\\
1. & $\begin{bmatrix}
0.83944 & 0.1465\\
0.071469 & 0.81019
\end{bmatrix}$ & 2.271  & -0.015311
\\
\\
\midrule
\\
2. &$\begin{bmatrix} 
0.72394 & 0.011739 \\ 
0.09617 & 0.82701
\end{bmatrix}$ & 2.2478  & -0.0084443
\\
\\
\midrule
\\
3. &$\begin{bmatrix}  
0.94433 & 0.0010596 \\ 
0.29449 & 0.92457
\end{bmatrix}$ & 2.2426 & -0.021666
\\
\\
\bottomrule
\end{tabular}
\caption{Numerical values of $\xi, \gamma, \text{ and } J_\gamma(\xi, \phi_k) \text{ when } J_\gamma(\xi, \phi_k) < 0$}
\end{table}
\clearpage

\item For \begin{equation*}
    \phi_3^0(x,y) = \begin{bmatrix}
    \displaystyle\frac{1}{100}\sin(\pi x), \displaystyle\frac{1}{100}\sin(\frac{3\pi}{2}y)
\end{bmatrix}
\end{equation*}

\begin{figure}[htp]\label{fig:smallermeshXiFun3}
    \centering
    \includegraphics[width=12cm]{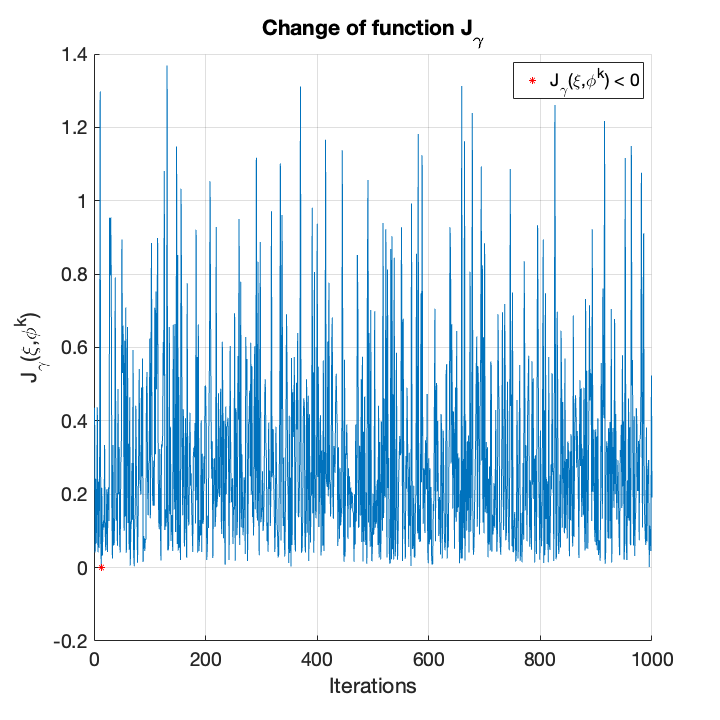}
    \caption{Change of function during the steepest descent on $\phi^0_3(x,y) = [
    \displaystyle\frac{1}{100}\sin(\pi x), \displaystyle\frac{1}{100}\sin(\displaystyle\frac{3\pi}{2}y)
    ]$ with changing $\xi$ and $\gamma$}
\end{figure}
\clearpage

\begin{table}\label{Table 6}
\centering
\begin{tabular}{llll}
\toprule
    & $\qquad\qquad\xi$ \qquad\qquad\qquad\qquad\qquad & $\quad\gamma$ \qquad\qquad\qquad\qquad& $J_\gamma(\xi, \phi_k)$  \\
\midrule
\\
1. & $\begin{bmatrix}
0.86869 & 0.088879\\
0.17486 & 0.88434
\end{bmatrix}$ & 2.3057 & -0.00058231
\\
\\
\bottomrule
\end{tabular}
\caption{Numerical values of $\xi, \gamma, \text{ and } J_\gamma(\xi, \phi_k) \text{ when } J_\gamma(\xi, \phi_k) < 0$}
\end{table}
\clearpage
\end{enumerate}
\clearpage
\subsubsection{$h = \frac{1}{10}$, $\xi$ with all its entries normal distributed between $0$ and $10$}
\begin{enumerate}
    \item For \begin{equation*}
        \phi_1^0(x,y) = \begin{bmatrix}
    \sin(x(x-1)y(y-1)), \sin(x(x-1)y(y-1))^2
\end{bmatrix}
    \end{equation*}
\begin{figure}[htp]\label{fig:mesh10XiFun1}
    \centering
    \includegraphics[width=12cm]{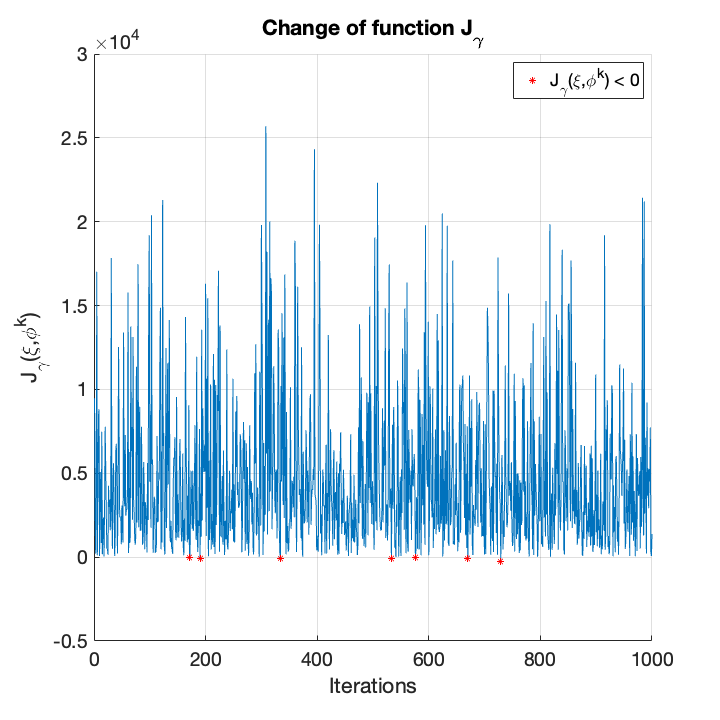}
    \caption{Change of function during the steepest descent on $\phi^0_1(x,y) = [
    \sin(x(x-1)y(y-1)), \sin(x(x-1)y(y-1))^2
    ]$ with changing $\xi$ and $\gamma$}
\end{figure}
\newpage

\begin{table}\label{Table 7}
\centering
\begin{tabular}{llll}
\toprule
    & $\qquad\qquad\xi$ \qquad\qquad\qquad\qquad\qquad & $\quad\gamma$ \qquad\qquad\qquad\qquad& $J_\gamma(\xi, \phi_k)$  \\
\midrule
\\
1. & $\begin{bmatrix}
7.1776 & 0.60596\\
2.9995 & 8.575
\end{bmatrix}$ & 2.2565  & -15.0015
\\
\\
\midrule
\\
2. &$\begin{bmatrix} 
5.9909 & 1.7873 \\ 
0.20139 & 8.4131
\end{bmatrix}$ & 2.2509 & -54.9543
\\
\\
\midrule
\\
3. &$\begin{bmatrix}  
6.783 & 1.3577 \\ 
0.60081 & 9.0508
\end{bmatrix}$ & 2.2064 & -87.8799
\\
\\
\midrule
\\
4. &$\begin{bmatrix}  
7.245 & 0.56314 \\ 
0.010342 & 5.7141
\end{bmatrix}$ & 2.1448 & -63.0402
\\
\\
\midrule
\\
5. &$\begin{bmatrix}  
9.0147 & 1.335 \\ 
0.62929 & 6.8561
\end{bmatrix}$ & 2.1315 & -2.2415
\\
\\
\midrule
\\
6. &$\begin{bmatrix}  
8.9732 & 0.35831 \\ 
1.0778 & 7.6424
\end{bmatrix}$ & 2.1024 & -102.372
\\
\\
\midrule
\\
7. &$\begin{bmatrix}  
9.7037 & 0.54265 \\ 
0.47065 & 9.7007
\end{bmatrix}$ & 2.0842 & -272.5885
\\
\\
\bottomrule
\end{tabular}
\caption{Numerical values of $\xi, \gamma, \text{ and } J_\gamma(\xi, \phi_k) \text{ when } J_\gamma(\xi, \phi_k) < 0$}
\end{table}
\clearpage

\item For \begin{equation*}
    \phi_2^0(x,y) = \begin{bmatrix}
    (x(x-1)y(y-1)), 0
\end{bmatrix}
\end{equation*}

\begin{figure}[htp]\label{fig:mesh10XiFun2}
    \centering
    \includegraphics[width=12cm]{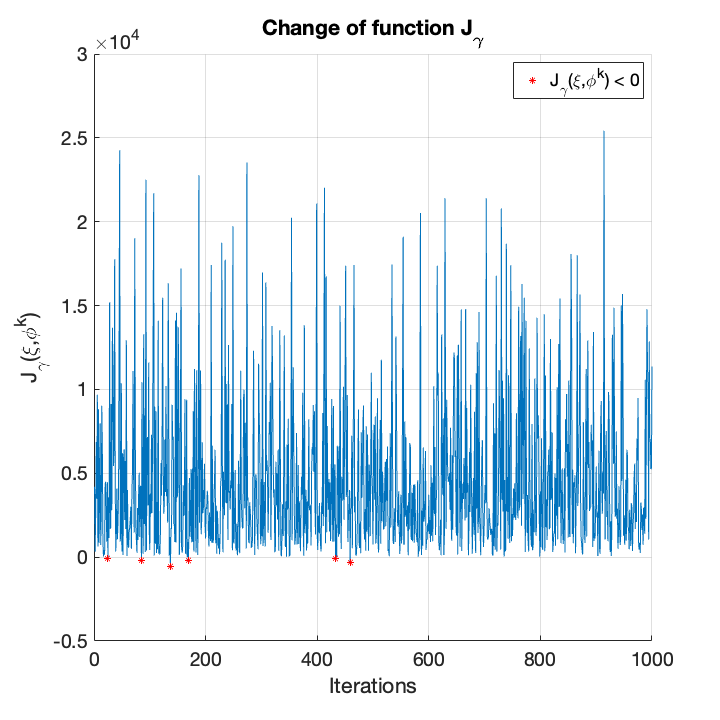}
    \caption{Change of function during the steepest descent on $\phi^0_2(x,y) = [
    (x(x-1)y(y-1)), 0
    ]$ with changing $\xi$ and $\gamma$}
\end{figure}
\newpage

\begin{table}\label{Table 8}
\centering
\begin{tabular}{llll}
\toprule
    & $\qquad\qquad\xi$ \qquad\qquad\qquad\qquad\qquad & $\quad\gamma$ \qquad\qquad\qquad\qquad& $J_\gamma(\xi, \phi_k)$  \\
\midrule
\\
1. & $\begin{bmatrix}
4.6987 &1.4696\\
0.014401 & 5.2997
\end{bmatrix}$ & 2.302  &  -53.0445
\\
\\
\midrule
\\
2. &$\begin{bmatrix}  
9.4925 & 1.5355 \\ 
0.25502 & 6.5724
\end{bmatrix}$ & 2.2834 &-181.2429
\\
\\
\midrule
\\
3. &$\begin{bmatrix}  
8.7277 & 0.97433 \\ 
0.65593 & 8.8613
\end{bmatrix}$ & 2.2673 &-570.7369
\\
\\
\midrule
\\
4. &$\begin{bmatrix}  
9.4153 & 1.616 \\ 
2.0106 & 8.3604
\end{bmatrix}$ & 2.2574 & -198.4102
\\
\\
\midrule
\\
5. &$\begin{bmatrix}  
5.9704 & 0.6489 \\ 
0.34235 & 5.4928
\end{bmatrix}$ & 2.1757 & -63.3355
\\
\\
\midrule
\\
6. &$\begin{bmatrix}  
8.4074 & 0.53192 \\ 
1.1414 & 9.1345
\end{bmatrix}$ &  2.1677 &-304.8817
\\
\\
\bottomrule
\end{tabular}
\caption{Numerical values of $\xi, \gamma, \text{ and } J_\gamma(\xi, \phi_k) \text{ when } J_\gamma(\xi, \phi_k) < 0$}
\end{table}
\clearpage

\item For \begin{equation*}
    \phi_3^0(x,y) = \begin{bmatrix}
    \displaystyle\frac{1}{100}\sin(\pi x), \displaystyle\frac{1}{100}\sin(\frac{3\pi}{2}y)
\end{bmatrix}
\end{equation*}

\begin{figure}[htp]\label{fig:mesh10XiFun3}
    \centering
    \includegraphics[width=12cm]{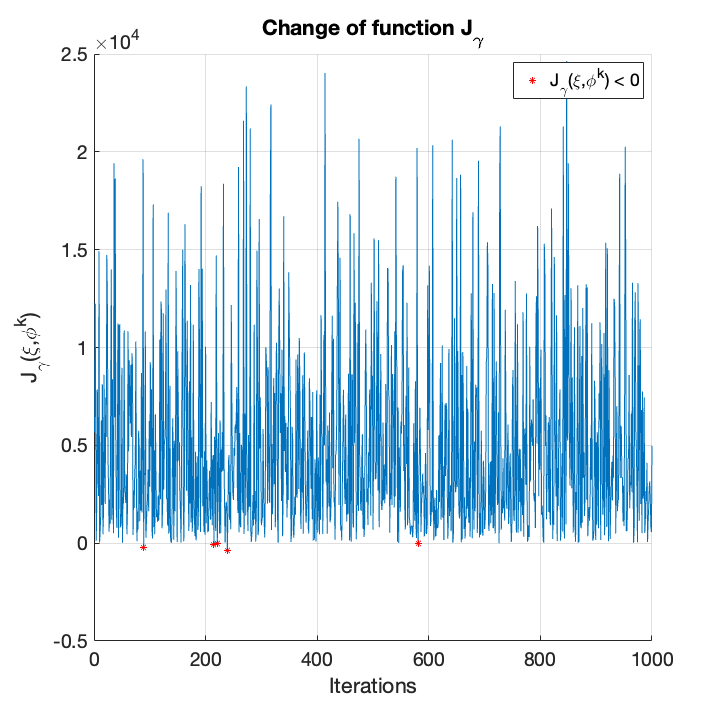}
    \caption{Change of function during the steepest descent on $\phi^0_3(x,y) = [
    \displaystyle\frac{1}{100}\sin(\pi x), \displaystyle\frac{1}{100}\sin(\displaystyle\frac{3\pi}{2}y)
    ]$ with changing $\xi$ and $\gamma$}
\end{figure}
\newpage
\begin{table}\label{Table 9}
\centering
\begin{tabular}{llll}
\toprule
    & $\qquad\qquad\xi$ \qquad\qquad\qquad\qquad\qquad & $\quad\gamma$ \qquad\qquad\qquad\qquad& $J_\gamma(\xi, \phi_k)$  \\
\midrule
\\
1. & $\begin{bmatrix}
6.5364 & 0.31994\\
1.4218 & 8.597
\end{bmatrix}$ & 2.2822 & -205.0655
\\
\\
\midrule
\\
2. &$\begin{bmatrix}  
7.2798 & 0.46976 \\ 
2.1372 & 6.3297
\end{bmatrix}$ & 2.2432 & -58.5336
\\
\\
\midrule
\\
3. &$\begin{bmatrix}  
6.9509 & 2.2428 \\ 
0.74368 & 6.6601
\end{bmatrix}$ & 2.2413 & -28.61
\\
\\
\midrule
\\
4. &$\begin{bmatrix}  
9.2135 & 2.6252 \\ 
0.11778 & 8.665
\end{bmatrix}$ & 2.2358 & -366.9857
\\
\\
\midrule
\\
5. &$\begin{bmatrix}  
8.7022 & 1.0632 \\ 
0.94471 & 7.2276
\end{bmatrix}$ & 2.1299 & -40.9998
\\
\\
\bottomrule
\end{tabular}
\caption{Numerical values of $\xi, \gamma, \text{ and } J_\gamma(\xi, \phi_k) \text{ when } J_\gamma(\xi, \phi_k) < 0$}
\end{table}
\clearpage
\end{enumerate}

\subsubsection{$h = \frac{1}{20}$, $\xi$ with all its entries normal distributed between $0$ and $10$}
\begin{enumerate}
    \item For \begin{equation*}
        \phi_1^0(x,y) = \begin{bmatrix}
    \sin(x(x-1)y(y-1)), \sin(x(x-1)y(y-1))^2
\end{bmatrix}
    \end{equation*}
\begin{figure}[htp]\label{fig:smallermesh10XiFun1}
    \centering
    \includegraphics[width=12cm]{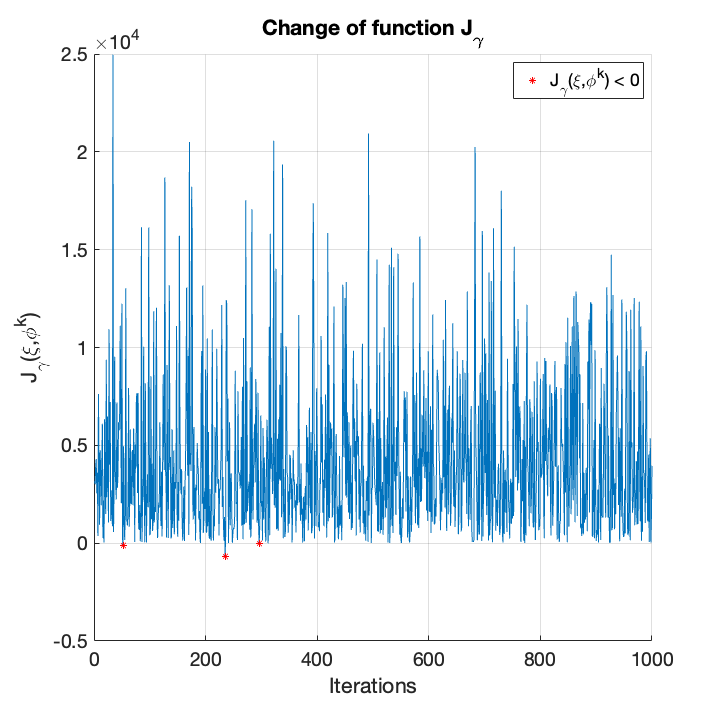}
    \caption{Change of function during the steepest descent on $\phi^0_1(x,y) = [
    \sin(x(x-1)y(y-1)), \sin(x(x-1)y(y-1))^2
    ]$ with changing $\xi$ and $\gamma$}
\end{figure}
\newpage

\begin{table}\label{Table 10}
\centering
\begin{tabular}{llll}
\toprule
    & $\qquad\qquad\xi$ \qquad\qquad\qquad\qquad\qquad & $\quad\gamma$ \qquad\qquad\qquad\qquad& $J_\gamma(\xi, \phi_k)$  \\
\midrule
\\
1. & $\begin{bmatrix}
8.3364 & 2.7385\\
1.4706 & 9.0459
\end{bmatrix}$ & 2.2936  & -101.9413
\\
\\
\midrule
\\
2. &$\begin{bmatrix} 
9.7842 & 0.51061 \\ 
1.6791 & 9.5567
\end{bmatrix}$ & 2.237 & -680.8128
\\
\\
\midrule
\\
3. &$\begin{bmatrix}  
5.5244 & 2.1469 \\ 
0.25663 & 6.2953
\end{bmatrix}$ & 2.2181 & -18.6345
\\
\\
\bottomrule
\end{tabular}
\caption{Numerical values of $\xi, \gamma, \text{ and } J_\gamma(\xi, \phi_k) \text{ when } J_\gamma(\xi, \phi_k) < 0$}
\end{table}
\clearpage

\item For \begin{equation*}
    \phi_2^0(x,y) = \begin{bmatrix}
    (x(x-1)y(y-1)), 0
\end{bmatrix}
\end{equation*}

\begin{figure}[htp]\label{fig:smallermesh10XiFun2}
    \centering
    \includegraphics[width=12cm]{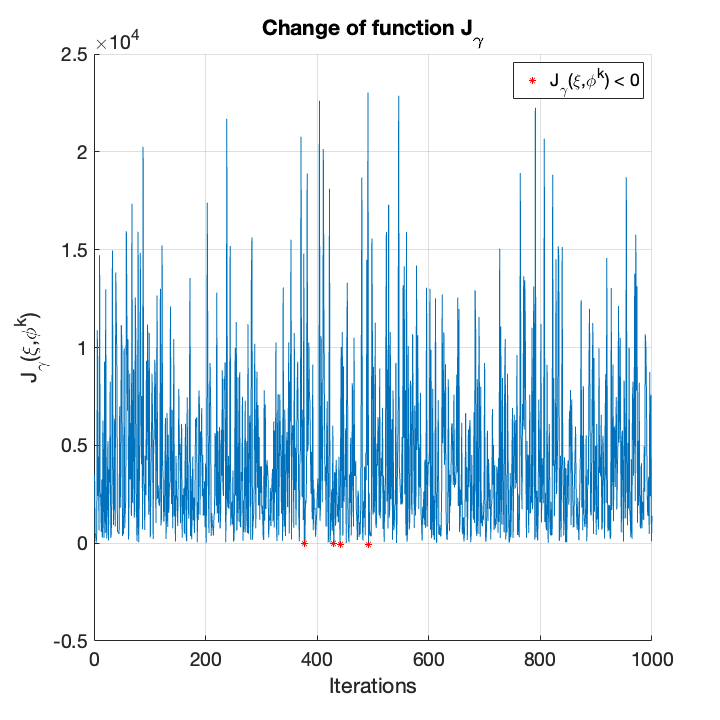}
    \caption{Change of function during the steepest descent on $\phi^0_2(x,y) = [
    (x(x-1)y(y-1)), 0
    ]$ with changing $\xi$ and $\gamma$}
\end{figure}
\newpage

\begin{table}\label{Table 11}
\centering
\begin{tabular}{llll}
\toprule
    & $\qquad\qquad\xi$ \qquad\qquad\qquad\qquad\qquad & $\quad\gamma$ \qquad\qquad\qquad\qquad& $J_\gamma(\xi, \phi_k)$  \\
\midrule
\\
1. & $\begin{bmatrix}
7.7158 &0.88689\\
1.4781 &5.9489
\end{bmatrix}$ & 2.1931 &  -9.7677
\\
\\
\midrule
\\
2. &$\begin{bmatrix}  
3.8016 & 0.8981 \\ 
0.063279 & 2.8256
\end{bmatrix}$ & 2.1767 &-0.57306
\\
\\
\midrule
\\
3. &$\begin{bmatrix}  
9.9608 & 1.4626 \\ 
1.1029 & 7.602
\end{bmatrix}$ & 2.1733 &-46.4173
\\
\\
\midrule
\\
4. &$\begin{bmatrix}  
7.3999& 1.4032 \\ 
0.39645 & 9.7827
\end{bmatrix}$ & 2.1575 & -67.3851
\\
\\
\bottomrule
\end{tabular}
\caption{Numerical values of $\xi, \gamma, \text{ and } J_\gamma(\xi, \phi_k) \text{ when } J_\gamma(\xi, \phi_k) < 0$}
\end{table}
\clearpage

\item For \begin{equation*}
    \phi_3^0(x,y) = \begin{bmatrix}
    \displaystyle\frac{1}{100}\sin(\pi x), \displaystyle\frac{1}{100}\sin(\frac{3\pi}{2}y)
\end{bmatrix}
\end{equation*}

\begin{figure}[htp]\label{fig:smallermesh10XiFun3}
    \centering
    \includegraphics[width=12cm]{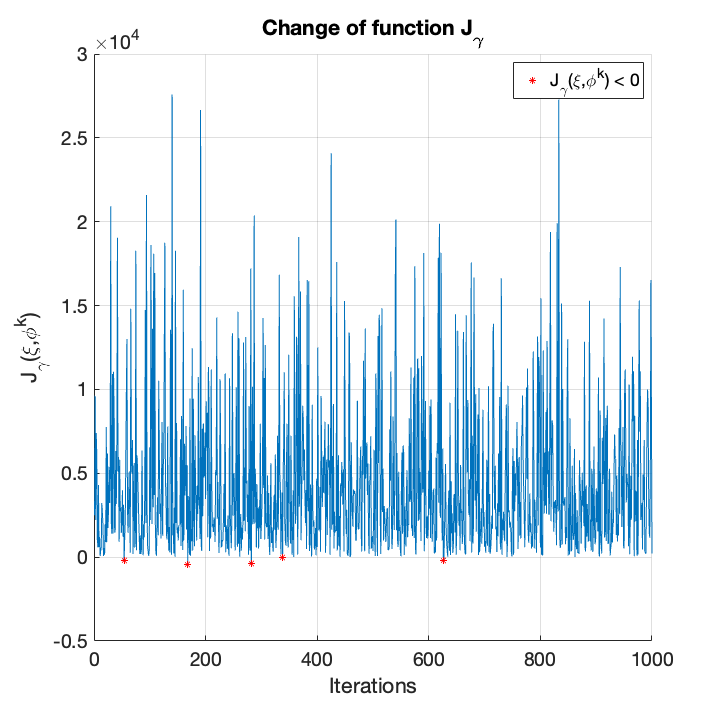}
    \caption{Change of function during the steepest descent on $\phi^0_3(x,y) = [
    \displaystyle\frac{1}{100}\sin(\pi x), \displaystyle\frac{1}{100}\sin(\displaystyle\frac{3\pi}{2}y)
    ]$ with changing $\xi$ and $\gamma$}
\end{figure}
\newpage

\begin{table}\label{Table 12}
\centering
\begin{tabular}{llll}
\toprule
    & $\qquad\qquad\xi$ \qquad\qquad\qquad\qquad\qquad & $\quad\gamma$ \qquad\qquad\qquad\qquad& $J_\gamma(\xi, \phi_k)$  \\
\midrule
\\
1. & $\begin{bmatrix}
8.6493 & 3.0288\\
1.0635 & 9.3453
\end{bmatrix}$ & 2.293 &-181.5976
\\
\\
\midrule
\\
2. &$\begin{bmatrix}  
8.3502 & 0.7933 \\ 
1.1585 & 8.4219
\end{bmatrix}$ & 2.2577 & -409.6368
\\
\\
\midrule
\\
3. &$\begin{bmatrix}  
9.0082 & 1.8819 \\ 
0.98139 & 9.5437
\end{bmatrix}$ & 2.2225 & -365.0452
\\
\\
\midrule
\\
4. &$\begin{bmatrix}  
5.8729 & 1.522 \\ 
0.21613 & 4.8101
\end{bmatrix}$ & 2.2051 & -16.003
\\
\\
\midrule
\\
5. &$\begin{bmatrix}  
8.6712 & 0.1889 \\ 
0.092423& 7.6643
\end{bmatrix}$ & 2.1157 & -179.8743
\\
\\
\bottomrule
\end{tabular}
\caption{Numerical values of $\xi, \gamma, \text{ and } J_\gamma(\xi, \phi_k) \text{ when } J_\gamma(\xi, \phi_k) < 0$}
\end{table}
\clearpage
\end{enumerate}

\section{Discussion}
\subsection{Discussion of Our Results}
The results we found are of several types and demonstrate that Morrey Conjecture is valid, at least numerically speaking. \hyperref[100sqrt3]{Figure 2} - \hyperref[fig:smallermesh10XiFun3]{Figure 15} represent the changes of Jensen integral $J_\gamma(\xi,\phi)$ as the iterations proceed for various values of $\xi \in \mathbb{R}^{2 \times 2}$ and $\gamma$. The oscillations observed in \hyperref[fig:randEpsFun1]{Figure 5} - \hyperref[fig:smallermesh10XiFun3]{Figure 15} do not characterize any default of the algorithm we designed and used; they appear because we used random  values of $\xi$: at each of the steepest descent algorithm iterations, a new value of the matrix $\xi$ is used. Doing so allows us to increase our chance of finding a matrix $\xi$, a value of $\gamma$ and a numerical mapping $\phi$ for which the Jensen inequality is violated. And indeed, for instance, for the values of $\xi = \begin{bmatrix}
0.79982& 0.099004\\
0.20266& 0.83511
\end{bmatrix}$ and $\gamma = 2.2723$, our numerical simulations show that the expression $J_\gamma(\xi,\phi)$ in Eq.(\ref{eqn:Jensen Inequality}) becomes negative, which violates the Jensen's inequality, and thereby confirming that the function is not quasi-convex. The rest of the analogous values of $\xi$ and $\gamma$ can be found from \hyperref[Table 1]{Table 1} - \hyperref[Table 12]{Table 12}. The three figures \hyperref[fig:randEpsFun1]{Figure 4} - \hyperref[fig:randEpsFun4]{Figure 6} differ from each other by the initial trial functions used to enter the steepest descent algorithm. We obtained the violation of the Jensen's inequality for initial guess mappings that are
\begin{equation}
  \phi^0 (x,y) =
    \begin{cases}
      & \phi^0_1(x,y) = \begin{bmatrix}
    \sin(x(x-1)y(y-1)), \quad \sin(x(x-1)y(y-1))^2
    \end{bmatrix}\\\
      & \phi^0_2(x,y) = \begin{bmatrix}
    (x(x-1)y(y-1)), \quad 0
    \end{bmatrix} \\
    & \phi^0_3(x,y) = \begin{bmatrix}
    \displaystyle\frac{1}{100}\sin(\pi x), \quad\displaystyle\frac{1}{100}\sin(\displaystyle\frac{3\pi}{2}y)
    \end{bmatrix}
    \end{cases}       
\end{equation}
The results would still hold if we use higher order of magnitude of the entries of $\xi$ and different mesh sizes as shown in \hyperref[fig:scaledXi]{Figure 8} where the entries of $\xi$ are multiplied by $20$ and in \hyperref[fig:smallermesh]{Figure 9} where we use a smaller mesh size of $\displaystyle\frac{1}{30}$.
$\\$

For fixed values of $\xi$, the values of the Jensen function at each iteration shown in \hyperref[100sqrt3]{Figure 2} and \hyperref[0000 figures]{Figure 3} show a decrease of the Jensen function until a minimal value after which it becomes a plateau or starts growing slowly. The oscillations observed in some of those figures maybe the results of numerical approximations that probably come from the discretization of the domain we used as it is the case in the numerical simulation of problem using finite element or finite difference methods. The decrease of the Jensen function is in the line of the steepest descent approach we used for the minimization problem we are solving. However, it is unclear the origin of the slow growth of this function observed in some of the figures after reaching a minimal point. It is also unclear why there is a sharp descent followed by a substantial growth.
\subsection{Comparison with Dacorogna \cite{dacorogna1990} and Gremaud's \cite{gremaud} Results}
Some comparisons of our results with Dacorogna \cite{dacorogna1990} and Gremaud \cite{gremaud} are at hand here since these two authors have worked with the same example function $f_\gamma$. In Gremaud's work \cite{gremaud}, a simulated annealing method was used but on the same example function $f_\gamma$. His simulation results indicate that rank-one convexity implies quasi-convexity for $\gamma = 1.1571;1.1553; 1.1550;1.1547$ (note that $\gamma$ in Gremaud's work has been scaled to half of the $\gamma$ we are using in our simulations). Gremaud's \cite{gremaud} results are aligned with the conclusions of the numerical simulations of Dacorogna et.al \cite{dacorogna1990}.
$\\$

In Dacorogna et.al. \cite{dacorogna1990}, a steepest gradient descent algorithm was used but with a crude approximation of the calculation of the gradient of $J_\gamma(\xi,\phi)$, and the use of linear interpolation functions to represent the functions $\phi$ the minimization problem relied on. Our algorithm and its modified version differ from that of Dacorogna et.al. \cite{dacorogna1990} by the fact that we calculated an exact analytic expression for the gradient of $J_\gamma(\xi,\phi)$, and we used a robust update of the step size in the steepest descent iteration. Dacorogna et.al. \cite{dacorogna1990} simulation results indicate that for $\xi = [1,0;0,\sqrt{3}]$ and four random functions $\phi$, $\gamma = 2.33$ corresponds to $J_\gamma = -44.05;-354.01;-440.55;-73.552$ which tells us that the function $f_\gamma$ is not quasi-convex but it is not rank-one convex (it would have been the case if $2 < \gamma <\displaystyle\frac{4}{\sqrt{3}} = 2.30940108$). Dacorogna et.al. \cite{dacorogna1990} reported also that for four different random values of $\xi$ and $\phi$, they got $J_\gamma(\xi,\phi) = 0$ up to the order of $10^{-12}$ when $\gamma = 2.31;2.32;2.3225;2.3275$ which also indicates that the simulation is not in favor of the Morrey Conjecture. Note that these values of $\gamma$ closely approach $\displaystyle\frac{4}{\sqrt{3}}$. For values of $\gamma$
that are very close to $\displaystyle\frac{4}{\sqrt{3}}$, our simulations, in contrast, reveal that $f_\gamma$ is non quasi-convex; for instance $\gamma=2.3017$. Even with $2<\gamma<2.3017$, we still values for $\gamma$ for which our simulations indicate that $f_\gamma$ is not quasi-convex. We believe that the reasons of this success are three-fold: (i) an accurate expression of the gradient of $J_\gamma(\xi,\phi)$; (ii) our initial guess functions $\phi$ of sinusoidal type; (iii) randomization of $\xi$ at each iteration and not only for four different random values of $\xi$ as in Dacorogna et.al \cite{dacorogna1990}. Our results were obtained for the entries of $\xi$ that are randomly distributed between $0$ and $1$ and $10 \times \xi$, and for meshsizes that are $h = \displaystyle\frac{1}{10}$ and $h = \displaystyle\frac{1}{20}$.
\section{Conclusion}
We report here the results of some numerical simulations we performed to address the Morrey Conjecture \cite{morrey1952quasi} problem. Unlike in the method proposed by Dacorogna \cite{dacorogna1990} and Gremaud \cite{gremaud}, where the gradient descent algorithm was either ``ad hoc’’ or unrelated to the minimization problem to be addressed, our approach used an exact vector gradient of the functional  $I(\phi) $ to be minimized over a Sobolev space which we defined previously. We derived an exact expression of this gradient. Then, we solved the minimization problem numerically by using the approximated values of the trial functions at each of the nodes on the mesh we used. Once the updated values of the trial functions at the nodes are obtained, we used them to check the Jensen's inequality which is associated to the quasi-convexity property of the function $f_\gamma$. The initial trial functions to enter the gradient descent iterative algorithm are chosen as oscillating functions for some fixed values of the $2 \times 2$ matrix $\xi$. By randomizing the entries of $\xi$, we successfully used for each of the iterations a new value of $\xi$.
With this modified algorithm, we obtain (i) a $\gamma$ for which Dacorogna function \cite{dacorogna1988counterexample} $f_\gamma$ is rank-one convex, (ii) some $\xi \in \mathbb{R}^{2\times 2}$, and a numerical $\phi^k$ for which the Jensen inequality is violated (the function is not quasi-convex), and thereby validating Morrey Conjecture.
$\\$
\section*{Software Availability}
A version of the code developed for this work is available at: \href{https://github.com/xdong99/Numerical-Quasiconvexity}{https://github.com/xdong99/Numerical-Quasiconvexity}.

\section*{Acknowledgements}
The authors would like to thank Dr. Wilfrid Gangbo for introducing them to the subject and for his useful comments on the first version of the paper. 
\clearpage
\bibliographystyle{unsrt}  
\bibliography{references}

\clearpage
\section*{Appendix: Calculation of the Gradient $\delta J_\gamma(\xi,\phi^k)$}\label{app: calculation grad}
Recall that $\phi\in W_0^{1,\infty}(\Omega;\mathbb{R}^2)$, which is a set of Lipshitz functions that vanish on $\partial\Omega$, we write
\begin{align}
    \phi &\coloneqq 
    \begin{bmatrix}
        \phi_1(x_1,x_2)\\
        \phi_2(x_1,x_2)
    \end{bmatrix}
\end{align}

It is also worth marking out notations that
\begin{align}
    \xi &= 
    \begin{bmatrix}
        \xi_{11} &\xi_{12}\\
        \xi_{21} &\xi_{22}
    \end{bmatrix},\\
    J_\gamma(\phi) &=  \int_\Omega f_\gamma(\xi + \nabla\phi)dx.
\end{align}

The steepest descent algorithm we use is in the form:
\begin{equation}
    \phi^{k+1} = \phi^k - \tau\delta J_\gamma(\phi^k)
\end{equation}
where $\tau$, the step size, is determined in each iteration. 
\\
\\
What's left to do is solving for $\delta J_\gamma(\phi^k).$ We notice the difficulty and complexity in finding this gradient of the functional $J_\gamma$ directly, so we proceed from writing out the gateaux derivative of $J_\gamma$ first.
\\
\\
The Gateaux differential $dJ(\phi;\psi)$ of $J_\gamma$ at $\phi \in \Omega$ in the direction of $\psi \in \mathbb{R}^2$ can be approached by
\begin{align}
\begin{split}
    J_\gamma(\phi + \varepsilon\psi) 
    &= \int_\Omega f_\gamma(\xi + \nabla\phi + \varepsilon\nabla\psi) - f_\gamma(\xi + \nabla\phi)dx \\
    &= \varepsilon\int_\Omega\sum_{i,j=1}^2\frac{\partial}{\partial \xi_{ij}}(\xi+\nabla\phi)\frac{\partial\psi_i}{\partial x_j}dx + o(\varepsilon) + J_\gamma(\phi).
\end{split}
\end{align}

Subtracting $J_\gamma(\phi)$ from both sides, the above expression, as we send $\epsilon$ to 0, yields $dJ(\phi;\psi)$
\begin{align}\label{dJphi}
\begin{split}
    dJ(\phi;\psi) &= \lim_{\varepsilon \to 0}\frac{J_\gamma(\phi + \varepsilon \psi) - J_\gamma(\phi)}{\varepsilon}\\
    &= \sum_{i,j=1}^2\int_\Omega\frac{\partial f_\gamma}{\partial\xi_{ij}}(\xi + \nabla\phi)\frac{\partial\psi_i}{\partial x_j}dx \\
    &= -\sum_{i,j=1}^2\int_\Omega\frac{\partial}{\partial x_j}\frac{\partial}{\partial \xi_{ij}}f_\gamma(\xi + \nabla\phi)\psi_i\ dx
\end{split}
\end{align}
which, by the definition of gateaux derivative, is also
\begin{equation}\label{gateaux}
    dJ(\phi;\psi) = \langle\delta J_\gamma(\phi),\psi\rangle \coloneqq \int_\Omega \delta J_\gamma(\phi(x))\psi(x)dx.
\end{equation}

From Eq.(\ref{dJphi}) and Eq.(\ref{gateaux}), it follows easily that
\begin{equation}\label{deltaJ}
    \delta J_\gamma(\phi) = -\sum_{i,j=1}^2\frac{\partial}{\partial x_j}\frac{\partial}{\partial \xi_{ij}}f_\gamma(\xi + \nabla\phi).
\end{equation}

Using the result in Eq.(\ref{deltaJ}), our algorithm can now be expressed as
\begin{equation}
    \phi^{k+1} = \phi^k + \tau\sum_{i,j=1}^2\frac{\partial}{\partial x_j}\frac{\partial}{\partial \xi_{ij}}f_\gamma(\xi + \nabla\phi^k).
\end{equation}

Since $\phi \in \mathbb{R}^2$, we perform gradient descent separately on each of the components
\begin{equation}
    \begin{cases}
    \vspace{0.25cm}
    \phi_1^{k+1} &= \phi_1^k + \tau\sum_{j=1}^2\frac{\partial}{\partial x_j}\frac{\partial}{\partial \xi_{1j}}f_\gamma(\xi + \nabla\phi^k)\\
    \vspace{0.25cm}
    &= \phi_1^k + \tau(\frac{\partial}{\partial x_1}\frac{\partial}{\partial \xi_{11}}f_\gamma(\xi + \nabla\phi^k) + \frac{\partial}{\partial x_2}\frac{\partial}{\partial \xi_{12}}f_\gamma(\xi + \nabla\phi^k))\\
    \vspace{0.25cm}
    \phi_2^{k+1} &= \phi_2^k + \tau\sum_{j=1}^2\frac{\partial}{\partial x_j}\frac{\partial}{\partial \xi_{1j}}f_\gamma(\xi + \nabla\phi^k)\\
    &= \phi_2^k + \tau(\frac{\partial}{\partial x_1}\frac{\partial}{\partial \xi_{12}}f_\gamma(\xi + \nabla\phi^k) + \frac{\partial}{\partial x_2}\frac{\partial}{\partial \xi_{22}}f_\gamma(\xi + \nabla\phi^k))\\
    \end{cases}
\end{equation}

There are four unknown partial derivatives in Eq.(\ref{update algorithm matlab}). We will show the complete steps for the first gradient. Since all four calculations are largely identical, we will only present the results of the other three.
\\
\\
Let $\nabla_j\phi_i$ denote $\displaystyle\frac{\partial}{\partial x_j}\phi_i$, then
\begin{align}
\begin{split}
    \xi + \nabla\phi &= 
        \begin{bmatrix}
            \xi_{11}+\nabla_1\phi_1 &\xi_{12}+\nabla_2\phi_1\\
            \xi_{21}+\nabla_1\phi_2 &\xi_{22}+\nabla_2\phi_2
        \end{bmatrix}
\end{split}
\end{align}

As $\xi + \nabla\phi$ is a regular $2\times 2$ matrix, it is clear that its determinant is
\begin{align}
\begin{split}
    det(\xi + \nabla\phi) = &(\xi_{11}+\nabla_1\phi_1)(\xi_{22}+\nabla_2\phi_2)\\
    &-(\xi_{12}+\nabla_2\phi_1)(\xi_{21}+\nabla_1\phi_2)
\end{split}
\end{align}
and its norm is
\begin{align}
\begin{split}
        \|\xi + \nabla\phi\|^2 = &(\xi_{11}+\nabla_1\phi_1)^2 + (\xi_{12}+\nabla_2\phi_1)^2 \\
          &+ (\xi_{21}+\nabla_1\phi_2)^2 + (\xi_{22}+\nabla_2\phi_2)^.
\end{split}
\end{align}

Let's start by computing
\begin{align}
\begin{split}
        \frac{\partial}{\partial x_1}\frac{\partial}{\partial \xi_{11}}f_\gamma(\xi + \nabla\phi).
\end{split}
\end{align}

To this end, we first consider the inner partial derivative
\begin{align}
\begin{split}
    \frac{\partial}{\partial \xi_{11}}f_\gamma(\xi + \nabla\phi)
\end{split}
\end{align}
which we calculate to be
\begin{align}
\begin{split}
    \frac{\partial}{\partial \xi_{11}}f_\gamma(\xi + \nabla\phi) &= \frac{\partial}{\partial \xi_{11}}\|\xi+\nabla\phi\|^4 - \frac{\partial}{\partial \xi_{11}}(\gamma\ \|\xi+\nabla\phi\|^2\ det(\xi+\nabla\phi)).
\end{split}
\end{align}

To make things cleaner, we calculate one term at a time,
\begin{align}\label{partial11ofnorm4}
\begin{split}
    \frac{\partial}{\partial \xi_{11}}\|\xi+\nabla\phi\|^4 &= \frac{\partial}{\partial \xi_{11}}(\|\xi + \nabla\phi\|^2)^2\\
    &= 2 \|\xi + \nabla\phi\|^2 \frac{\partial}{\partial \xi_{11}}(\|\xi + \nabla\phi\|^2)\\
\end{split}
\end{align}
where
\begin{align}\label{partial11ofnorm2}
\begin{split}
    \frac{\partial}{\partial \xi_{11}}(\|\xi + \nabla\phi\|^2) &= 2(\xi_{11}+\nabla_1\phi_1).
\end{split}
\end{align}

From Eq.(\ref{partial11ofnorm4}) and Eq.(\ref{partial11ofnorm2}),
\begin{align}\label{partial11ofnorm4 final}
\begin{split}
    \frac{\partial}{\partial \xi_{11}}\|\xi+\nabla\phi\|^4 &= 2 \|\xi + \nabla\phi\|^2(2(\xi_{11}+\nabla_1\phi_1))\\
    &= 4(\xi_{11}+\nabla_1\phi_1)\|\xi + \nabla\phi\|^2.
\end{split}
\end{align}

Now by the chain rule, we have
\begin{align}\label{chain rule for the other one}
\begin{split}
    \frac{\partial}{\partial \xi_{11}}(\gamma\ \|\xi+\nabla\phi\|^2\ det(\xi+\nabla\phi)) 
    &= \gamma \frac{\partial}{\partial \xi_{11}}(\|\xi + \nabla\phi\|^2\ det(\xi+\nabla\phi))\\
    &= \gamma\frac{\partial}{\partial \xi_{11}}(\|\xi + \nabla\phi\|^2)det(\xi + \nabla\phi) \\
    &\ \ \ \ + \gamma \|\xi + \nabla\phi\|^2 \frac{\partial}{\partial \xi_{11}}det(\xi + \nabla\phi)
\end{split}
\end{align}
where
\begin{align}\label{chain rule 1}
\begin{split}
    \gamma\frac{\partial}{\partial \xi_{11}}(\|\xi + \nabla\phi\|^2)det(\xi + \nabla\phi)
    &= 2\ \gamma\ (\xi_{11}+\nabla_1\phi_1)\ det(\xi+\nabla\phi)
\end{split}
\end{align}
and
\begin{align}\label{chain rule 2}
\begin{split}
    \gamma \|\xi + \nabla\phi\|^2 \frac{\partial}{\partial \xi_{11}}det(\xi + \nabla\phi)
    &= \gamma\ (\xi_{22}+\nabla_2\phi_2)\ \|\xi + \nabla\phi\|^2.
\end{split}
\end{align}

Combining the two equations Eq.(\ref{chain rule 1}) and Eq.(\ref{chain rule 2}), we can obtain
\begin{align}\label{chain rule for the other one final}
\begin{split}
    \frac{\partial}{\partial \xi_{11}}(\gamma\ \|\xi+\nabla\phi\|^2\ det(\xi+\nabla\phi)) 
    &= 2\ \gamma\ (\xi_{11}+\nabla_1\phi_1)\ det(\xi+\nabla\phi) \\
    &\quad + \gamma\ (\xi_{22}+\nabla_2\phi_2)\ \|\xi + \nabla\phi\|^2.
\end{split}
\end{align}

Now using the results from Eq.(\ref{partial11ofnorm4 final}) and Eq.(\ref{chain rule for the other one final}), we have
\begin{align}\label{partialXi11of f}
\begin{split}
    \frac{\partial}{\partial \xi_{11}}f_\gamma(\xi + \nabla\phi) 
    = &4(\xi_{11}+\nabla_1\phi_1)\|\xi + \nabla\phi\|^2 - (2\ \gamma\ (\xi_{11}+\nabla_1\phi_1)\\
    &det(\xi+\nabla\phi)+\gamma\ (\xi_{22}+\nabla_2\phi_2)\ \|\xi + \nabla\phi\|^2).
\end{split}
\end{align}

The remaining task is to find the partial derivative of Eq.(\ref{partialXi11of f}) with respect to $x_1$:
\begin{equation}\label{intermidiate gradient}
    \begin{split}
        \frac{\partial}{\partial x_1}\frac{\partial}{\partial \xi_{11}}f_\gamma(\xi + \nabla\phi) 
        =& \frac{\partial}{\partial x_1}(4(\xi_{11}+\nabla_1\phi_1)\|\xi + \nabla\phi\|^2) \\
        &- (\frac{\partial}{\partial x_1}(2\ \gamma\ (\xi_{11}+\nabla_1\phi_1)det(\xi+\nabla\phi))\\
        &+ \frac{\partial}{\partial x_1}(\gamma\ (\xi_{22}+\nabla_2\phi_2)\ \|\xi + \nabla\phi\|^2)) \\
    \end{split}
\end{equation}

We separate the work as before. First, we have
\begin{equation}\label{5.27}
    \begin{split}
        \frac{\partial}{\partial x_1}(4(\xi_{11}+\nabla_1\phi_1)\|\xi + \nabla\phi\|^2)
        = &4(\frac{\partial}{\partial x_1}(\xi_{11}+\nabla_1\phi_1))\|\xi + \nabla\phi\|^2 \\
          &+ 4(\xi_{11}+\nabla_1\phi_1)(\frac{\partial}{\partial x_1}\|\xi + \nabla\phi\|^2)
    \end{split}
\end{equation}
where the first term
\begin{equation}\label{5.28}
    \begin{split}
        4(\frac{\partial}{\partial x_1}(\xi_{11}+\nabla_1\phi_1))\|\xi + \nabla\phi\|^2
        &= 4\ \nabla_{11}\phi_1\|\xi + \nabla\phi\|^2
    \end{split}
\end{equation}
and we leave the second term unexpanded for now.
\\
\\
For the second term in Eq.(\ref{intermidiate gradient}), we obtain
\begin{equation}\label{5.29}
    \begin{split}
        \frac{\partial}{\partial x_1}2\ \gamma\ (\xi_{11}+\nabla_1\phi_1)&det(\xi+\nabla\phi)\\
        =&2\ \gamma\ (\frac{\partial}{\partial x_1}(\xi_{11}+\nabla_1\phi_1))det(\xi+\nabla\phi)) \\
         &+2\ \gamma\ (\xi_{11}+\nabla_1\phi_1)(\frac{\partial}{\partial x_1}det(\xi+\nabla\phi)) \\
        =&2\ \gamma\ \nabla_{11}\phi_1\ det(\xi+\nabla\phi) \\
        &+2\ \gamma\ (\xi_{11}+\nabla_1\phi_1))(\frac{\partial}{\partial x_1}det(\xi+\nabla\phi))
    \end{split}
\end{equation}

Similarly, for the third term in Eq.(\ref{intermidiate gradient}), we derive
\begin{equation}\label{5.30}
    \begin{split}
        \frac{\partial}{\partial x_1}(\gamma\ (\xi_{22}+\nabla_2\phi_2)\ &\|\xi + \nabla\phi\|^2))\\
        =&\gamma\ (\frac{\partial}{\partial x_1}(\xi_{22}+\nabla_2\phi_2))\|\xi + \nabla\phi\|^2 \\
         &+ \gamma\ (\xi_{22}+\nabla_2\phi_2)(\frac{\partial}{\partial x_1}\|\xi + \nabla\phi\|^2)\\
        =&\gamma\ \nabla_{21}\phi_2\ \|\xi + \nabla\phi\|^2\\
         &+\gamma\ (\xi_{22}+\nabla_{2}\phi_2)(\frac{\partial}{\partial x_1}\|\xi + \nabla\phi\|^2)
    \end{split}
\end{equation}

Combining Eq.(\ref{5.27}), Eq.(\ref{5.28}), Eq.(\ref{5.29}), and Eq.(\ref{5.30}), we conclude that
\begin{equation}
    \begin{split}
        \frac{\partial}{\partial x_1}&\frac{\partial}{\partial \xi_{11}}f_\gamma(\xi + \nabla\phi)\\
        = &4\ \nabla_{11}\phi_1\|\xi + \nabla\phi\|^2 + 4(\xi_{11}+\nabla_{1}\phi_1)\frac{\partial}{\partial x_1}\|\xi + \nabla\phi\|^2 \\
          &- 2\ \gamma\ \nabla_{11}\phi_1\ det(\xi+\nabla\phi) - 2\ \gamma\ (\xi_{11}+\nabla_1\phi_1))(\frac{\partial}{\partial x_1}det(\xi+\nabla\phi)) \\
          &-\gamma\ \nabla_{21}\phi_2\ \|\xi + \nabla\phi\|^2 -\gamma\ (\xi_{22}+\nabla_{2}\phi_2)(\frac{\partial}{\partial x_1}\|\xi + \nabla\phi\|^2)
    \end{split}
\end{equation}

We present the results of the other three gradients as well:
\begin{align}
\begin{split}
 1)\ \frac{\partial}{\partial x_2}&\frac{\partial}{\partial \xi_{12}}f_\gamma(\xi + \nabla\phi)\\ 
    &= 4\ \nabla_{22}\phi_1\|\xi + \nabla\phi\|^2 +4(\xi_{12}+\nabla_{2}\phi_1)\frac{\partial}{\partial x_2}\|\xi + \nabla\phi\|^2 \\
    &- 2\ \gamma\ \nabla_{22}\phi_1 det(\xi+\nabla\phi) - 2\ \gamma\ (\xi_{12}+\nabla_2\phi_1))\frac{\partial}{\partial x_2}det(\xi+\nabla\phi)\\
    &- \gamma\ \nabla_{12}\phi_2 \|\xi + \nabla\phi\|^2 -\gamma\ (\xi_{21}+\nabla_{1}\phi_1)\frac{\partial}{\partial x_2}\|\xi + \nabla\phi\|^2   
\end{split}
\end{align}

\begin{align}
\begin{split}
 2)\ \frac{\partial}{\partial x_1}&\frac{\partial}{\partial \xi_{21}}f_\gamma(\xi + \nabla\phi)\\ 
    &= 4\ \nabla_{11}\phi_2\|\xi + \nabla\phi\|^2 +4(\xi_{21}+\nabla_{1}\phi_2)\frac{\partial}{\partial x_1}\|\xi + \nabla\phi\|^2 \\
    &- 2\ \gamma\ \nabla_{11}\phi_2 det(\xi+\nabla\phi) - 2\ \gamma\ (\xi_{21}+\nabla_1\phi_2))\frac{\partial}{\partial x_1}det(\xi+\nabla\phi)\\
    &- \gamma\ \nabla_{21}\phi_1\|\xi + \nabla\phi\|^2 -\gamma\ (\xi_{12}+\nabla_{2}\phi_1)\frac{\partial}{\partial x_1}\|\xi + \nabla\phi\|^2  
\end{split}
\end{align}

\begin{align}
\begin{split}
 3)\ \frac{\partial}{\partial x_2}&\frac{\partial}{\partial \xi_{22}}f_\gamma(\xi + \nabla\phi)\\
    &= 4\ \nabla_{22}\phi_2\|\xi + \nabla\phi\|^2 +4(\xi_{22}+\nabla_{2}\phi_2)\frac{\partial}{\partial x_2}\|\xi + \nabla\phi\|^2\\
    &- 2\ \gamma\ \nabla_{22}\phi_2 det(\xi+\nabla\phi) - 2\ \gamma\ (\xi_{22}+\nabla_2\phi_2))\frac{\partial}{\partial x_2}det(\xi+\nabla\phi)\\
    &- \gamma\ \nabla_{12}\phi_1\|\xi + \nabla\phi\|^2 -\gamma\ (\xi_{11}+\nabla_{1}\phi_1)\frac{\partial}{\partial x_2}\|\xi + \nabla\phi\|^2   
\end{split}
\end{align}

For the sake of readability, we keep some partial derivatives unexpanded in the above results. For the common terms that appear in all the gradients above, we can write separate functions to calculate them in MATLAB$\copyright$\cite{MATLAB:2021} to avoid redundancy. Their expansions are given as below:
\begin{align}
\begin{split}
        \frac{\partial}{\partial x_1}\|\xi + \nabla\phi\|^2
        &= 2((\xi_{11}+\nabla_1\phi_1)\nabla_{11}\phi_1 +(\xi_{12}+\nabla_2\phi_1)\nabla_{21}\phi_1 \\
        &\quad + (\xi_{21}+\nabla_1\phi_2)\nabla_{11}\phi_2 + (\xi_{22}+\nabla_2\phi_2)\nabla_{21}\phi_2)\\
        \frac{\partial}{\partial x_2}\|\xi + \nabla\phi\|^2
        &= 2((\xi_{11}+\nabla_1\phi_1)\nabla_{12}\phi_1 +(\xi_{12}+\nabla_2\phi_1)\nabla_{22}\phi_1 \\
        &\quad + (\xi_{21}+\nabla_1\phi_2)\nabla_{12}\phi_2 + (\xi_{22}+\nabla_2\phi_2)\nabla_{22}\phi_2)\\
        \frac{\partial}{\partial x_1}det(\xi + \nabla\phi)
        &= \nabla_{11}\phi_1(\xi_{22}+\nabla_{2}\phi_2) + \nabla_{21}\phi_2(\xi_{11}+\nabla_{1}\phi_1)\\ 
        &\quad -\nabla_{11}\phi_2(\xi_{12}+\nabla_{2}\phi_1) - \nabla_{21}\phi_1(\xi_{21}+\nabla_{1}\phi_2)\\
        \frac{\partial}{\partial x_2}det(\xi + \nabla\phi)
        &= \nabla_{12}\phi_1(\xi_{22}+\nabla_{2}\phi_2) + \nabla_{22}\phi_2(\xi_{11}+\nabla_{1}\phi_1)\\ 
        &\quad -\nabla_{12}\phi_2(\xi_{12}+\nabla_{2}\phi_1) - \nabla_{22}\phi_1(\xi_{21}+\nabla_{1}\phi_2)
\end{split}
\end{align}

\end{document}